\newif\ifpreprint
\preprinttrue  

\ifpreprint

\documentclass[11pt]{article}
\usepackage{fullpage}
\title{Verification of Sequential Convex Programming for Parametric Non-convex Optimization}

\author{Rajiv Sambharya, Nikolai Matni, George Pappas}
\date{%
    University of Pennsylvania\\
    \today
}

\usepackage[round,authoryear]{natbib}
\usepackage{amsmath,enumitem,mathtools,amsthm,amssymb}

\newtheorem{theorem}{Theorem}

\newtheorem{lemma}[theorem]{Lemma}
\newtheorem{example}[theorem]{Example}
\newtheorem{proposition}[theorem]{Proposition}

\newtheorem{assumption}[theorem]{Assumption}

\usepackage[colorlinks,
            linkcolor=cyan,
            citecolor=cyan,
            urlcolor=cyan,
            bookmarks]{hyperref}

\else 

\documentclass[final,onefignum,onetabnum]{siamart250211}

\usepackage{lipsum}
\usepackage{amsfonts}
\usepackage{graphicx}
\usepackage{epstopdf}
\usepackage{bibentry}
\ifpdf
  \DeclareGraphicsExtensions{.eps,.pdf,.png,.jpg}
\else
  \DeclareGraphicsExtensions{.eps}
\fi

\usepackage{enumitem}
\setlist[enumerate]{leftmargin=.5in}
\setlist[itemize]{leftmargin=.5in}

\newsiamremark{remark}{Remark}
\newsiamremark{hypothesis}{Hypothesis}
\crefname{hypothesis}{Hypothesis}{Hypotheses}
\newsiamthm{claim}{Claim}

\headers{Verification of Sequential Convex Programming}{R. Sambharya, N. Matni, and G. Pappas}

\title{Verification of Sequential Convex Programming for Parametric Non-convex Optimization\thanks{Submitted to the editors 11/26/2025.}}

\author{Rajiv Sambharya\thanks{University of Pennsylvania, Philadelphia, PA (\email{sambhar9@seas.upenn.edu}, \\ \email{nmatni@seas.upenn.edu}, \email{pappasg@seas.upenn.edu})} \and Nikolai Matni\footnotemark[2] \and George Pappas\footnotemark[2]} 

\usepackage{amsopn}

\let\citep\cite
\let\citet\cite

\usepackage{amssymb}
\usepackage{cite}
\newtheorem{assumption}[theorem]{Assumption}

\usepackage[skip=8pt]{caption}
\usepackage{graphicx}
\usepackage{subcaption}
\usepackage{hyperref}
\fi
\usepackage{tikz}
\usetikzlibrary{positioning}
\usepackage{verbatim}
\usepackage{multicol}
\usepackage{algorithm}
\usepackage[noend]{algpseudocode} 
\usepackage{bm}

\usepackage{subcaption}
\usepackage{makecell}
\usepackage{csvsimple}	
\usepackage{booktabs}   
\usepackage[flushleft]{threeparttable}  
\usepackage[round-mode=places, round-integer-to-decimal, round-precision=4,
    table-format = 1.4,
    table-number-alignment=center,
    round-integer-to-decimal]{siunitx}
\usepackage{adjustbox}  

\newcommand{\bnote}[1]{}
\renewcommand{\bnote}[1]{\textcolor{red}{\textbf{B. #1}}}

\usepackage{etoolbox}

\DeclareMathOperator*{\argmin}{argmin}

\newcommand*{\startlegend}{-0.2}
\newcommand*{\enlegend}{0.7}

\DeclareDocumentCommand{\T}{ O{z} O{} }{T\IfValueT{#2}{(#1,\theta_{#2})}\IfNoValueT{#2}{(#1,\theta)}}
\DeclareDocumentCommand{\Tj}{ O{k} O{z} O{} }{T^{#1}\IfValueT{#3}_{\theta_{#3}}{(#2)}}

\DeclareDocumentCommand{\CB}{ O{} }{C_{B_{\IfValueTF{#1}{\theta_{#1}}{\theta}}}}
\DeclareDocumentCommand{\RB}{ O{} }{R_{B_{\IfValueTF{#1}{\theta_{#1}}{\theta}}}}

\newcommand{\eg}{{\it e.g.}}
\newcommand{\ie}{{\it i.e.}}

\newcommand{\ones}{\mathbf 1}
\newcommand{\reals}{{\mbox{\bf R}}}

\newcommand{\symm}{{\mbox{\bf S}}}  

\newcommand{\Sec}{Section}
\newcommand{\Subsec}{Subsection}
\newcommand{\Thm}{Theorem}

\newcommand{\Eqn}{Equation}

\newcommand{\cblock}[3]{
  \hspace{-1.5mm}
  \begin{tikzpicture}
    [
    node/.style={square, minimum size=10mm, thick, line width=0pt},
    ]
    \node[fill={rgb,255:red,#1;green,#2;blue,#3}] () [] {};
  \end{tikzpicture}%
}

\usetikzlibrary{shapes.geometric}

\newcommand{\linecircle}[3]{%
  \begin{tikzpicture}[baseline={(0,-.1)}]
    \draw[draw=none, fill={rgb,255:red,#1;green,#2;blue,#3}](0.25, 0) circle (.08);
    \draw[line width=0.8pt, color={rgb,255:red,#1;green,#2;blue,#3}] (\startlegend,0) -- (\enlegend,0);
  \end{tikzpicture}%
  }

\newcommand{\linecirclehollow}[3]{%
  \begin{tikzpicture}[baseline={(0,-.1)}]
    \draw[line width=0.8pt, draw={rgb,255:red,#1;green,#2;blue,#3}](0.25, 0) circle (.08);
    \draw[line width=0.8pt, dotted, color={rgb,255:red,#1;green,#2;blue,#3}] (\startlegend,0) -- (\enlegend,0);
  \end{tikzpicture}%
  }

\newcommand{\linesquarehollow}[3]{%
  \begin{tikzpicture}[baseline={(0,-.1)}]
    \draw[line width=0.8pt, draw={rgb,255:red,#1;green,#2;blue,#3}](0.179, -.071) rectangle (.321,.071);
    \draw[line width=0.8pt, dotted, color={rgb,255:red,#1;green,#2;blue,#3}] (\startlegend,0) -- (\enlegend,0);
  \end{tikzpicture}%
  }

\newcommand{\linedotted}[3]{%
\begin{tikzpicture}[baseline={(0,-.1)}]
  \draw[line width=1.5pt, dotted, color={rgb,255:red,#1;green,#2;blue,#3}] (\startlegend,0) -- (\enlegend,0);
\end{tikzpicture}%
}

\newcommand{\linesquare}[3]{%
  \begin{tikzpicture}[baseline={(0,-.1)}]
    \draw[draw=none, fill={rgb,255:red,#1;green,#2;blue,#3}](0.179, -.071) rectangle (.321,.071);
    \draw[line width=0.8pt, color={rgb,255:red,#1;green,#2;blue,#3}] (\startlegend,0) -- (\enlegend,0);
  \end{tikzpicture}%
  }

\newcommand{\linedowntri}[3]{
  \begin{tikzpicture}[baseline={(0,-.1)}]
    \coordinate (A) at (0.25,-.09);
    \coordinate (B) at (.17,.06);
    \coordinate (C) at (.33,.06);
    \draw[draw=none, fill={rgb,255:red,#1;green,#2;blue,#3}] (A) -- (B) -- (C) -- cycle;
    \draw[line width=0.8pt, color={rgb,255:red,#1;green,#2;blue,#3}] (\startlegend,0) -- (\enlegend,0);
  \end{tikzpicture}%
}

\newcommand{\lineuptri}[3]{
  \begin{tikzpicture}[baseline={(0,-.1)}]
    \coordinate (A) at (0.25,.09);
    \coordinate (B) at (.17,-.06);
    \coordinate (C) at (.33,-.06);
    \draw[draw=none, fill={rgb,255:red,#1;green,#2;blue,#3}] (A) -- (B) -- (C) -- cycle;
    \draw[line width=0.8pt, color={rgb,255:red,#1;green,#2;blue,#3}] (\startlegend,0) -- (\enlegend,0);
  \end{tikzpicture}%
}

\newcommand{\lineuptrihollow}[3]{
  \begin{tikzpicture}[baseline={(0,-.1)}]
    \coordinate (A) at (0.25,.09);
    \coordinate (B) at (.17,-.06);
    \coordinate (C) at (.33,-.06);
    \draw[line width=0.8pt, draw={rgb,255:red,#1;green,#2;blue,#3}] (A) -- (B) -- (C) -- cycle;
    \draw[line width=0.8pt, dotted, color={rgb,255:red,#1;green,#2;blue,#3}] (\startlegend,0) -- (\enlegend,0);
  \end{tikzpicture}%
}

\newcommand{\linedowntrihollow}[3]{
  \begin{tikzpicture}[baseline={(0,-.1)}]
    \coordinate (A) at (0.25,-.09);
    \coordinate (B) at (.17,.06);
    \coordinate (C) at (.33,.06);
    \draw[line width=0.8pt, draw={rgb,255:red,#1;green,#2;blue,#3}] (A) -- (B) -- (C) -- cycle;
    \draw[line width=0.8pt, dotted, color={rgb,255:red,#1;green,#2;blue,#3}] (\startlegend,0) -- (\enlegend,0);
  \end{tikzpicture}%
}

\newenvironment{talign*}
 {\let\displaystyle\textstyle\csname align*\endcsname}
 {\endalign}

\newcommand{\tableheader}{
  \begin{tabular}{@{}c@{}}Tol. \\\end{tabular}
  &
\begin{tabular}{@{}c@{}}Nesterov\end{tabular}
&
\begin{tabular}{@{}c@{}}Nearest \\ Neighbor\end{tabular}
&
\begin{tabular}{@{}c@{}}L2WS \\ $N=10$\end{tabular}
&
\begin{tabular}{@{}c@{}}L2WS \\ $N=10000$\end{tabular}
&
\begin{tabular}{@{}c@{}}LM \\ $N=10$\end{tabular}
&
\begin{tabular}{@{}c@{}}LM \\ $N=10000$\end{tabular}
&
\begin{tabular}{@{}c@{}}LAH \end{tabular}
&
\begin{tabular}{@{}c@{}}LAH \\ Robust\end{tabular}
\\}

\newcommand{\colnames}{\colA & \colB & \colC &\colD &\colE & \colF & \colG & \colH & \colI}

\newcommand{\myCSVReaderRed}[1]{%
\tableheader
      \midrule
    \csvreader[
        head to column names,
        late after line=\\
    ]{#1}{ 
        accuracies=\colA,
        cold_start_red=\colC,
        nearest_neighbor_red=\colD,
        maml_red=\colO,
        obj_k0_red=\colE,
        obj_k5_red=\colF,
        obj_k15_red=\colG,
        obj_k30_red=\colH,
        obj_k60_red=\colI,
        reg_k0_red=\colJ,
        reg_k5_red=\colK,
        reg_k15_red=\colL,
        reg_k30_red=\colM,
        reg_k60_red=\colN,
    }{\colnames}
}

\newcommand{\myCSVReaderRedAlt}[2]{%
    \tableheaderAlt{#1}
      \midrule
    \ifstrequal{#1}{MAML}{%
    \csvreader[
        head to column names,
        late after line=\\
    ]{#2}{ 
        accuracies=\colA,
        cold_start_red=\colC,
        nearest_neighbor_red=\colD,
        maml_red=\colO,
        obj_k0_red=\colE,
        obj_k1_red=\colF,
        obj_k5_red=\colG,
        obj_k15_red=\colH,
        obj_k60_red=\colI,
        reg_k0_red=\colJ,
        reg_k1_red=\colK,
        reg_k5_red=\colL,
        reg_k15_red=\colM,
        reg_k60_red=\colN,
    }{\colnames}
    }
    {
      \csvreader[
        head to column names,
        late after line=\\
    ]{#2}{ 
        accuracies=\colA,
        cold_start_red=\colC,
        nearest_neighbor_red=\colD,
        prev_sol_red=\colO,
        obj_k0_red=\colE,
        obj_k1_red=\colF,
        obj_k5_red=\colG,
        obj_k15_red=\colH,
        obj_k60_red=\colI,
        reg_k0_red=\colJ,
        reg_k1_red=\colK,
        reg_k5_red=\colL,
        reg_k15_red=\colM,
        reg_k60_red=\colN,
    }{\colnamesprevsol}
    }
}

\newcommand{\useCSVReaderRed}[2]{%
    \ifdefined\currentCSVReaderRed
        \expandafter\csname\currentCSVReaderRed\endcsname{#1}{#2}%
    \else
        \PackageWarning{YourPackage}{No CSV reader command defined, defaulting to \myCSVReaderRed}%
        \myCSVReaderRed{#1}{#2}%
    \fi
}

\makeatletter
\renewcommand{\eqref}[1]{\textup{\tagform@{\ref{#1}}}}
\makeatother

\begin{document}
\maketitle

\begin{abstract}
  We introduce a verification framework to exactly verify the worst-case performance of sequential convex programming (SCP) algorithms for parametric non-convex optimization. 
  The verification problem is formulated as an optimization problem that maximizes a performance metric (\eg, the suboptimality after a given number of iterations) over parameters constrained to be in a parameter set and iterate sequences consistent with the SCP update rules.
  Our framework is general, extending the notion of SCP to include both conventional variants such as trust-region, convex-concave, and prox-linear methods, and algorithms that combine convex subproblems with rounding steps, as in relaxing and rounding schemes.
  Unlike existing analyses that may only provide local guarantees under limited conditions, our framework delivers global worst-case guarantees—quantifying how well an SCP algorithm performs across all problem instances in the specified family.
  Applications in control, signal processing, and operations research demonstrate that our framework provides, for the first time, global worst-case guarantees for SCP algorithms in the parametric setting.
\end{abstract}

\newcommand{\myparagraph}[1]{%
  \paragraph{#1\ifpreprint\unskip.\fi}%
}

\newcommand{\myparagraphstar}[1]{%
  \par\vspace{1ex}\noindent\textbf{#1\ifpreprint.\fi}\hspace{1em plus 0.5em}%
}


\section{Introduction}\label{sec:intro}
In this paper, we are interested in solving the parametric non-convex optimization problem
\begin{equation}\label{prob:main}
  \begin{array}{ll}
  \mbox{minimize} & f(z,x) \\
  \mbox{subject to}  & z \in \Omega(x),
  \end{array}
\end{equation}
where $z \in \reals^n$ is the decision variable, $x \in \mathcal{X} \subseteq \reals^d$ is the \emph{problem parameter}, $\mathcal{X}$ is a known set, $\Omega(x) \subseteq \reals^n$ is the possibly non-convex feasible region, and $f : \reals^n \times \reals^d \rightarrow \reals$ is the possibly non-convex objective with respect to $z$.
In many settings, we must repeatedly solve problem~\eqref{prob:main} with a varying problem parameter $x$.
For instance, in model predictive control~\citep{borrelli_mpc_book,mpc_primal_active_sets}, we repeatedly solve optimal control problems with varying initial states.
In energy systems, we repeatedly solve optimal power flow problems with varying loads and renewable generation~\citep{ml_opf,warm_start_power_flow}.
In signal processing, we repeatedly solve sparse coding problems with varying noisy measurements~\citep{lista,alista}.
In operations research, such parametric structure arises in optimization problems over networks (\eg, traffic assignment, routing, or graph-based flow problems) where the underlying topology remains fixed but demands or costs vary across instances~\citep{still2018lectures}.

Global approaches such as branch-and-bound~\citep{lawler1966branch,minlo} and branch-and-cut~\citep{padberg1991branch,stubbs1999branch} can, in principle, provide globally optimal solutions to problem~\eqref{prob:main}. 
However, these methods scale poorly with problem size and are often not practical in time-sensitive settings where real-time optimization is required~\citep{online_milliseconds,bengio_machine_2021}.
Instead, it is common to rely on faster, heuristic approaches such as \emph{sequential convex programming} \citep{gill2005snopt,dinh2010local,NoceWrig06,boyd2008sequential} in which a sequence of convex problems is solved.
Convex optimization enjoys both a rich theoretical foundation, with strong guarantees of global optimality~\citep{cvxbook}, and a mature ecosystem of efficient algorithms~\citep{boyd2011admm,prox_algos,mon_primer,lscomo} and solvers~\citep{goulart2024clarabel,osqp,scs_quadratic} that make it practical for applications where real-time optimization is required.

Many methods fall under this broad category of sequential convex programming. 
One widely used approach is the trust-region method~\citep{NoceWrig06,byrd2000trust}, where at each iteration a local convex model is optimized within a bounded region around the current iterate. 
Another example is the convex-concave procedure (CCP) for difference-of-convex (DC) programming problems~\citep{lipp2016variations,yuille2001concave,shen2016disciplined}, where the concave part of the objective or constraints is linearized to yield a convex subproblem. 
For composite optimization, the prox-linear method~\citep{drusvyatskiy2019efficiency,drusvyatskiy2017proximal} is a popular approach. At each iteration, it replaces the smooth inner map with its linearization and evaluates the convex (possibly non-smooth) outer function on this approximation, thereby reducing each step to a convex subproblem.
It is common to use variants of these approaches, for example, by using penalized procedures or hyperparameters (\eg, the trust-region size) that vary over the iterations~\citep{lipp2016variations}.

Relaxing and rounding schemes are another popular way to solve problems of the form~\eqref{prob:main} when the feasible set includes discrete constraints. 
A common example is to relax a non-convex constraint to obtain a convex problem and then round the resulting solution~\citep{goemans1995improved,ageev2004pipage}.
The work of~\citet{diamond2018general} systematizes this approach with the relax-round-polish method: first, a relaxation of the non-convex problem is solved, then the discrete variables are rounded, and finally the remaining continuous variables are re-optimized by solving a convex problem.
Like SCP, these methods rely on solving a sequence of convex optimization problems, but they differ in that they include explicit non-convex projection steps.

While both conventional SCP methods and related heuristic methods that combine convex subproblems with rounding steps are widely used in practice, they lack guarantees on their global worst-case performance~\citep{diamond2018general,boyd2008sequential,lipp2016variations}.
Such limitations can be unacceptable in systems where it is essential to \emph{verify} the global quality of the candidate solution of these algorithms over all admissible parameters in the set $\mathcal{X}$.
At the same time, existing analyses typically do not take advantage of the parametric structure of problem~\eqref{prob:main}. 
Exploiting this structure creates an opportunity to construct global guarantees in cases where none are otherwise available.

\subsection{Contributions}
In this paper, we develop a framework to verify the global worst-case performance of sequential convex programming methods for parametric non-convex optimization.
While SCP is traditionally understood to involve only convex subproblems, we adopt a broader definition that also encompasses algorithms combining convex subproblems with simple non-convex projection steps (\eg, in rounding schemes).
Our verification method is meant to be solved offline (\ie, in a time-insensitive setting), so that the algorithm can then be used online when the problem parameter is seen (\ie, in a time-sensitive setting) with certified worst-case guarantees.
In this work, we focus on the setting of non-convex quadratically-constrained quadratic programs with potentially additional discrete constraints (\ie, binary or sparsity constraints).
Our key contributions are as follows:
\begin{itemize}[left=5pt]
  \item {\bf Verification framework for global performance guarantees.} 
  We introduce a verification framework to certify the global worst-case performance of SCP methods for parametric non-convex optimization. 
  In this framework, verification is posed as an optimization problem that maximizes a performance metric over problem parameters constrained to a specified set and iterate sequences consistent with the SCP update rules for a given number of iterations.
  We focus on three complementary aspects of performance: (i) suboptimality with respect to the optimal value of problem~\eqref{prob:main} in cases where feasibility of the candidate solution can be guaranteed, (ii) the worst-case level of constraint violation of the candidate solution when feasibility cannot be guaranteed, and (iii) feasibility of downstream convex subproblems which contain only linear constraints.
  \item {\bf Applicability across sequential convex programming algorithms.} 
  Our framework can be applied across a range of SCP algorithms including conventional variants such as trust-region, convex-concave, and prox-linear methods, and algorithms that combine convex subproblems with rounding steps such as relax-round-polish.
  The key idea to our approach is to encode the algorithm's steps—both convex quadratic programs (QPs) and rounding operations—as constraints in the verification problem through their optimality conditions.
  In the most general case, this yields a verification problem that is a mixed-integer quadratically-constrained quadratic program that can be solved using modern solvers.
  \item {\bf Numerical examples.} Through a broad variety of numerical examples from control, signal processing, and operations research, we demonstrate the utility of our framework in certifying worst-case guarantees on the performance of SCP algorithms.
  Our approach is able to provide exact numerical guarantees where no known analytic methods are able to.
  Our results reveal a broad spectrum of algorithmic behaviors, from certifiable convergence to optimality to cases where the worst-case guarantee plateaus short of global optimality, and even complete failure to make progress toward finding a feasible solution beyond the first iteration.
  These examples also show how the framework can provide further insights to design algorithms, such as hyperparameter selection and initialization choice.
\end{itemize}

\myparagraph{Layout of the paper}
In \Sec~\ref{sec:related_work}, we review related work.
In \Sec~\ref{sec:scp_algos}, we introduce several SCP algorithms including the trust-region method, the CCP, the prox-linear method, and the relax-round-polish method.
In \Sec~\ref{sec:verification}, we present our verification framework to obtain global worst-case performance guarantees.
In \Sec~\ref{sec:numerical_experiments}, we showcase the utility of this verification framework with many numerical examples.
Finally, in \Sec~\ref{sec:conclusion}, we conclude.

\myparagraph{Notation}
We use $\reals^n$ to denote the space of $n$-dimensional real vectors, $\reals_+^n$ the set of $n$-dimensional with non-negative entries, and $\reals_{++}^n$ the set of $n$-dimensional with positive real numbers.
We use $\symm^n$ to denote the space of $n \times n$-symmetric matrices, $\symm_+^n$ the set of $n \times n$-positive semidefinite symmetric matrices, and $\symm_{++}^n$ the set of $n \times n$-positive definite symmetric matrices.
We denote the all ones vector of length $m$ with $\mathbf{1_m}$, the all zeros vector of length $m$ with $\mathbf{0_m}$, and the $n \times n$ identity matrix with $I_n$.
For a set $S \subseteq \reals^n$ and point $v \in \reals^n$, we define the indicator function $\mathcal{I}_S(v) = 0$ if $v \in S$ and $\infty$ otherwise. 
For a vector $v \in \reals^n$, we let $v_+ = \max(v, 0)$ element-wise.

\section{Related work}\label{sec:related_work}

\myparagraph{Conventional sequential convex programming}
Sequential convex programming \citep{gill2005snopt,dinh2010local,NoceWrig06,boyd2008sequential,byrd2000trust,bonalli2019gusto} conventionally refers to an algorithmic approach that solves non-convex problems by repeatedly constructing and solving convex approximations. 
Classical examples include the trust-region method~\citep{NoceWrig06,byrd2000trust,boyd2008sequential}, the CCP for DC programming~\citep{tao1997convex,lipp2016variations,yuille2001concave}, and the prox-linear method for composite minimization~\citep{drusvyatskiy2019efficiency,lewis2016proximal}. 
In certain cases, some of these methods come with local guarantees~\citep{drusvyatskiy2019efficiency,bonalli2019gusto}, but in all cases, constructing tight global worst-case guarantees remains a challenge.
Our framework can exactly certify the worst-case performance of such SCP methods in the parametric setting of~\eqref{prob:main}.


\myparagraph{Relaxation and rounding methods}
Relaxation and rounding schemes are widely used to tackle non-convex problems by replacing difficult discrete constraints with relaxed convex constraints, followed by either deterministic~\citep{ageev2004pipage,nannicini2012rounding} or randomized~\citep{goemans1995improved,raghavan1987randomized} rounding to produce candidate solutions. 
A systematic formulation of this approach for deterministic rounding is the relax-round-polish framework~\citep{diamond2018general}, which combines a convex relaxation, a non-convex rounding step, and a subsequent polishing step that re-optimizes the continuous, convex variables.
Such methods often yield strong empirical performance, but they lack worst-case guarantees on suboptimality and feasibility~\citep{diamond2018general}. 
While SCP is traditionally understood to involve only convex subproblems, we adopt a broader definition that also encompasses algorithms combining convex subproblems with rounding steps (\eg, in relax-round-polish).
Our framework also provides a way to exactly certify the worst-case performance of these methods.

\myparagraph{Guarantees in parametric optimization}
Recently, several works have studied tight performance bounds for algorithms that solve parametric optimization problems.
The works of~\citet{ranjan2024exact,perfverifyqp} are most similar to ours: they verify the worst-case performance of first-order methods for parametric convex optimization problems by encoding the algorithm steps as constraints in the verification problem.
In contrast, our work verifies the performance of SCP algorithms for parametric non-convex optimization problems.
Another key distinction is that in our setting, the algorithmic steps are defined \emph{implicitly} through convex subproblems rather than through explicit update rules. 
Moreover, while global worst-case guarantees exist in convex optimization, no such guarantees are available in our non-convex setting.

Other methods provide \emph{probabilistic guarantees} on the performance of algorithms to solve parametric optimization problems~\citep{sambharya2025data,huang2025data}; however, these methods rely on the assumption that the problem parameters are independently and identically distributed (i.i.d.).
This assumption is often hard to verify and may not be realistic in many scenarios (\eg, when problems are solved sequentially as in control applications~\citep{borrelli_mpc_book}).
Moreover, these guarantees only hold with high-probability---meaning that performance can still be arbitrarily poor even if it occurs rarely. 


\myparagraph{Finite-iteration guarantees in learning to optimize}
An adjacent research area, learning to optimize (L2O)~\citep{l2o,amos_tutorial,bengio_machine_2021}, focuses on learning to accelerate algorithms for parametric optimization.
L2O methods have been applied to a range of non-convex settings—for example, learning branching rules or the discrete variables for mixed-integer optimization~\citep{bengio_machine_2021,khalil_learning_2016,online_milliseconds,cauligi2021coco}, and learning warm starts for SCP~\citep{banerjee2020learning}—where they often deliver impressive empirical results.
A central challenge in L2O, however, is guaranteeing performance after a finite number of iterations of an algorithm~\citep{l2o,amos_tutorial}.
One line of work addresses this gap with probabilistic generalization bounds, showing that a learned optimizer likely performs well on unseen problem instances~\citep{learn_algo_steps,Sucker2024LearningtoOptimizeWP,l2ws,saravanos2024deep,sambharya2025data,balcan_gen_guarantees}. 
Another line develops worst-case, finite-iteration guarantees, but these results hold only for convex problems~\citep{sambharya2025accel,banert2021accelerated,martin2025learning}.
Although our work also considers the parametric setting, it differs fundamentally in that we do not learn an optimizer.
Instead, we exactly verify the worst-case finite-iteration performance of SCP methods for parametric non-convex optimization—providing guarantees in a regime where existing L2O worst-case bounds do not apply.

\myparagraph{Neural network verification}
Our approach is closely related in spirit to neural network verification, which seeks to certify that for every admissible input (\eg, an image within a bounded perturbation set), the network's output satisfies a desired property (\eg, correct classification)~\citep{tjeng2017evaluating,ceccon2022omlt,anderson2020strong}.
The key similarity between the two lies in formulating verification as an optimization problem that searches for the most adversarial instance while encoding the computation steps (either the forward pass of a neural network or an SCP algorithm) as constraints.

\section{Sequential convex programming algorithms}\label{sec:scp_algos}
In this section, we introduce several SCP algorithms that fit within our verification framework.
We focus on the trust-region method in \Subsec~\ref{subsec:trust_scp}, the penalized CCP in \Subsec~\ref{subsec:pen_convex_concave}, the prox-linear method in \Subsec~\ref{subsec:prox_linear}, and the relax-round-polish method in \Subsec~\ref{subsec:relax_round_polish}.
In each type of algorithm, we first write the non-convex optimization problem under consideration (which is a special case of problem~\eqref{prob:main}), and then present the algorithm steps.
All methods can be described by the iterations $z^{k+1} \in s^k(z^k, x)$, where $z^k$ is the $k$-th iterate and $s^k : \reals^n \times \reals^d \rightrightarrows \reals^n$ is a (possibly multi-valued) operator.
In this paper, we focus on the quite general case where the non-convex problems are mixed-integer quadratically-constrained quadratic programs.
For the algorithms we consider\footnote{except for the relax-round-polish method which has no initial point}, the choice of initialization $z^0$ and hyperparameters can significantly influence the quality of the candidate solution produced by the algorithm.

\subsection{Trust-region method}\label{subsec:trust_scp}
In this subsection, we consider the trust-region method to solve the problem
\begin{equation}\label{prob:scp}
  \begin{array}{ll}
  \mbox{minimize} & f(z,x) \\
  \mbox{subject to} & g(z,x) \leq 0, \quad  h(z,x) = 0,
  \end{array}
\end{equation}
where $z \in \reals^n$ is the decision variable, $x \in \reals^d$ is the problem parameter, the functions $f : \reals^n \times \reals^d \rightarrow \reals$ and $g :\reals^n \times \reals^d \rightarrow \reals^{m}$ may be non-convex in $z$, and $h : \reals^n \times \reals^d \rightarrow \reals^{p}$ may be non-affine in $z$. 
The trust-region method iteratively builds local convex approximations around the current iterate and solves the resulting subproblems subject to a trust-region constraint~\citep{boyd2008sequential}.
Formally, given a current iterate $z^k$, we construct the trust-region subproblem
\begin{equation}\label{prob:tr}
  z^{k+1} \in s(z^k, x)
  = \begin{aligned}[t]
      &\text{argmin}   && \hat{f}(z, x, z^k) \\
      &\text{subject to} && \hat{g}(z, x, z^k) \le 0, \quad \hat{h}(z, x, z^k) = 0 \\
      &                  && \|z - z^k\|_\infty \le \rho_k,
    \end{aligned}
\end{equation}
where $\hat{f}(z, x, z^k)$ and $\hat{g}(z, x, z^k)$ are convex approximations of $f$ and $g$ around $z^k$, $\hat{h}(z, x, z^k)$ is an affine approximation of $h$ around $z^k$, and $\rho_k > 0$ is the trust-region radius at the $k$-th iteration. 
The trust-region constraint prevents the next iterate from moving too far from the current iterate, ensuring that the convexified problem is a meaningful local model of the non-convex problem.

\subsection{Penalized convex-concave procedure}\label{subsec:pen_convex_concave}
In this subsection, we consider the penalized CCP to solve the DC program
\begin{equation}\label{prob:dc}
  \begin{array}{ll}
  \mbox{minimize} & f_0(z,x) - g_0(z,x) \\
  \mbox{subject to} & f_i(z,x) - g_i(z,x) \leq 0, \quad i=1,\dots,m,
  \end{array}
\end{equation}
where $z \in \reals^n$ is the decision variable and $f_i : \reals^n \times \reals^d \rightarrow \reals$ and $g_i : \reals^n \times \reals^d \rightarrow \reals$ are convex functions in $z$ for $i=0,\dots,m$.
Problem~\eqref{prob:dc} is convex only if the $g_i$ functions are \emph{affine} in $z$.
The penalized CCP involves solving a series of optimization subproblems of the form
\begin{equation}\label{prob:pen_ccp}
  z^{k+1} \in s(z^k, x)
  = \begin{aligned}[t]
      &\text{argmin}   && f_0(z, x) - \hat{g}_0(z, x, z^k) + \tau_k \mathbf{1}_m^T s \\
      &\text{subject to} && f_i(z, x) - \hat{g}_i(z, x, z^k) \le s_i, \quad i = 1, \dots, m, \\
      &                  && s \ge 0,
    \end{aligned}
\end{equation}
where $z \in \reals^n$ (the original decision variable) and $s \in \reals^m$ (introduced slack variables) are the decision variables, and $\tau_k > 0$ is a hyperparameter that weights the constraint violations at the $k$-th iteration.
At each iterate $z^k$, the non-affine functions are replaced with their local affine approximations $\hat{g}_i(z,x,z^k) = g_i(z,x) + \nabla_z g_i(z^k, x)^T (z - z^k)$ for $i=0,\dots,m$.
By introducing slacks and penalizing their violation in the objective, the penalized CCP ensures that each subproblem remains feasible while still driving the iterates toward satisfaction of the original non-convex constraints.

\subsection{Prox-linear method}\label{subsec:prox_linear}
In this subsection, we focus on the prox-linear method, an algorithm that is specialized to minimize composite functions of the form
\begin{equation}\label{prob:prox_linear}
  \begin{array}{ll}
  \mbox{minimize} & g(z,x) + h(c(z,x), x),
  \end{array}
\end{equation}
where $g: \reals^n \times \reals^d \rightarrow \reals \cup \{+\infty\}$ is convex in $z$, $h: \reals^p \times \reals^d \rightarrow \reals$ is convex in $z$, and $c: \reals^n \times \reals^d \rightarrow \reals^p$ is smooth with respect to $z$.
The prox-linear method~\citep{drusvyatskiy2017proximal,drusvyatskiy2019efficiency} consists of solving the following convex subproblem at each iteration:
\begin{equation}\label{prob:prox_linear_subproblem}
  z^{k+1} \in s(z^k, x) = \begin{array}{ll}
  \mbox{argmin} & g(z,x) + h(c(z^k,x) + \nabla c(z^k,x)^T (z - z^k)) + \frac{\rho}{2} \|z - z^k\|_2^2,
  \end{array}
\end{equation}
where $z$ is the decision variable and $\rho > 0$ is a regularization hyperparameter.
Problem~\eqref{prob:prox_linear_subproblem} is convex in $z$ since convex functions composed with affine mappings are convex \citep[\Subsec~3.3.2]{cvxbook}.
We focus on the case where $h$ is the absolute value function and $c$ is a non-convex quadratic function---a case which has many applications~\citep{drusvyatskiy2017proximal}.
Under certain assumptions, the iterates of the prox-linear method are known to converge to a point that is nearly stationary~\citep{drusvyatskiy2019efficiency}.

\subsection{Relax-round-polish method}\label{subsec:relax_round_polish}
We now turn to an algorithm that is not typically viewed as an SCP method, due 
to the presence of a non-convex projection step: the relax-round-polish 
algorithm~\citep{diamond2018general}.  
We apply it to the non-convex problem
\begin{equation}\label{prob:rrp}
  \begin{array}{ll}
  \mbox{minimize} & (1/2)w^T P_w w + c_w^T w + (1/2)v^T P_v v + c_v^T v \\
  \mbox{subject to} & A w + B v \leq d, \quad v \in \mathcal{Z}, \\
  \end{array}
\end{equation}
where the decision variable is $z=(w,v)$ with $w \in \reals^{n_w}$ and 
$v \in \reals^{n_v}$.
The problem parameter is given by all of the problem data: \ie, $x = (P_w, c_w, P_v, c_v, A, B, d)$.
Moreover, $P_w \succeq 0$ and $P_v \succeq 0$, so problem~\eqref{prob:rrp} is a convex problem except for the non-convex constraint $v \in \mathcal{Z}$.
In this paper, we focus on two non-convex sets:
\begin{itemize}[left=5pt]
    \item Binary constraints: $\mathcal{Z} = \{0,1\}^{n_v}$.
    \item Sparsity constraints: $\mathcal{Z} = \{v \in \reals^{n_v} \mid \|v\|_0 \leq k\}$ for a given integer $k$.
\end{itemize}
The relax-round-polish method consists of the following three steps.
\begin{itemize}[left=5pt]
  \item {\bf Relax}: Solve a convex relaxation of problem~\eqref{prob:rrp}.
  For binary constraints, we replace $v \in \{0,1\}^{n_v}$ with $v \in [0,1]^{n_v}$.
  For sparsity constraints, we soften the condition by solving problem~\eqref{prob:rrp} without the non-convex constraint but with an added penalty of $\lambda \|v\|_1$ for a chosen hyperparameter $\lambda > 0$.
  \item {\bf Round}: Take the solution of the relaxed problem and project the variable $v$ to the non-convex set.
  For binary variables, this means rounding each entry of $v$ to either $0$ or $1$.
  For sparsity constraints, this means keeping the $k$ entries of largest magnitude and setting the other entries to $0$.
  \item {\bf Polish}: With the discrete variables $v$ fixed, solve the resulting convex subproblem in the decision variable $w$.
  This step sharpens the continuous variables $w$ while enforcing consistency with the discrete variables $v$.
\end{itemize}
Each of these three steps can be written in the form of the iterations $z^{k+1} \in s^k(z^k,x)$, where each operator $s^k$ is different.
In this case of relax-round-polish, there is no need for an initialization $z^0$.
While the relax and round steps are feasible optimization problems, there is no guarantee that the polish convex subproblem is feasible~\citep{diamond2018general}.

\section{Verification framework}\label{sec:verification}

In this section, we introduce our verification framework for certifying the global worst-case performance of SCP methods for parametric non-convex optimization.
In \Subsec~\ref{sec:formulating_verification_problem}, we formulate the verification problem as an optimization problem that maximizes a performance metric subject to constraints that i) encode the SCP steps, ii) enforce the parameter is in a set $\mathcal{X}$, iii) enforce that the initial point is in a set $S$, and iv) impose optimality on the true solution to problem~\eqref{prob:main}.
In \Subsec~\ref{subsec:algorithm_steps}, we show how to encode the algorithm steps (the convex QPs and rounding steps) as constraints in the verification problem.
Then in the next three subsections, we focus on three specific cases: verifying suboptimality when feasibility of the final iterate can be guaranteed in \Subsec~\ref{subsec:nonconvex_optimality}, verifying the level of constraint violation when feasibility cannot be guaranteed in \Subsec~\ref{subsec:constraint_violation}, and verifying the feasibility of convex subproblems with linear constraints in \Subsec~\ref{subsec:feas}.
Finally, in \Subsec~\ref{subsec:inexact}, we show how to incorporate inexact solves to the convex subproblems into our framework.

\subsection{Formulating the verification problem}\label{sec:formulating_verification_problem}
Before presenting the verification formulations, we detail the algorithmic updates, initialization sets, and performance metrics that together specify how the algorithms are run and evaluated.

\myparagraph{Algorithm steps}
All of the algorithms outlined in \Sec~\ref{sec:scp_algos} can be represented by the iterations
\begin{equation*}
  z^{k+1} \in s^k(z^k, x), \quad k=0,\dots,K-1,
\end{equation*}
where $s^k : \reals^n \times \reals^d \rightrightarrows \reals^n$ is a (possibly multi-valued) operator and $K$ is the number of iterations.
The number of steps $K$ is either dictated by i) the prescribed number of steps within an algorithm (\eg, exactly three in relax-round-polish) or ii) the number of iterations one has time to run (\eg, in the case of trust-region and prox-linear methods).
The operators $s^k$ need not be identical across iterations: different steps in an algorithm may correspond to solving different convex subproblems or applying different rounding rules.
Moreover, we model each $s^k$ as a multi-valued operator since all of the algorithm steps we consider may have more than one solutions.
Throughout, we assume that each subproblem admits at least one solution; issues of subproblem infeasibility are treated separately later in \Subsec~\ref{subsec:feas}.

\myparagraph{Initial points}
Most of the algorithms we consider require an initial point $z^0$ which can significantly affect the quality of the resulting candidate solution~\citep{boyd2008sequential}.
We model the initialization by assuming that $z^0$ belongs to some known non-empty set $S$.
For cold-started algorithms, the set of initial points is $S = \{{\bf 0}_n\}$.
When a heuristic warm start $z^{\rm ws} \in \reals^n$ is available, we simply take $S = \{z^{\rm ws}\}$.
For algorithms whose first step does not depend on $z^0$ (\eg, relax-round-polish), the specification of $S$ is omitted.

\myparagraph{Performance metrics}
We are interested in analyzing the worst-case performance in terms of some performance metric $\phi(z^K,z^\star,x)$ where $z^K$ is a candidate solution and $z^\star$ is the optimal solution to the original non-convex problem.
We focus on three key aspects of performance in our verification framework: suboptimality, the level constraint violation, and subproblem feasibility.
\begin{itemize}[left=5pt]
  \item {\bf Suboptimality.} To provide guarantees on suboptimality, we let 
  \begin{equation}\label{eq:suboptimality}
  \phi(z^K, z^\star, x) = f(z^K,x) - f(z^\star,x),
  \end{equation}
  where $f(z,x)$ is the objective function of the original non-convex problem~\eqref{prob:main}, and $z^\star$ is the solution parametrized by parameter $x$.
  In this case, it is necessary that the final iterate $z^K$ be feasible, which is the case in many scenarios (see \Subsec~\ref{subsec:nonconvex_optimality} for details).
  \item {\bf Level of constraint violation.} Given constraints of the form $\Omega(x) = \{z \mid g(z,x) \leq 0\}$, the square of the $\ell_2$-norm of violations is quantified as 
   \begin{equation}\label{eq:feas_viol}
    \phi(z^K, z^\star, x) = \sum_i (\max \{g_i(z^K,x), 0\})^2.
  \end{equation}
  We could easily adapt this metric to other norms, such as the $\ell_1$-norm or the $\ell_\infty$-norm of violations.
  Moreover, the metric~\eqref{eq:feas_viol} could be adapted to handle equality constraints.
  \item {\bf Subproblem feasibility.} 
  Sometimes, we cannot guarantee that the convex subproblems are feasible.
  For example, in the relax-round-polish method, once the discrete variables are fixed, there is no guarantee that there exists continuous variables that can satisfy the constraints.
  We show how to use the Farkas Lemma in \Subsec~\ref{subsec:feas} to check feasibility of downstream convex QPs.
\end{itemize}

\myparagraph{The verification problem}
Using the algorithm steps, initial points, and performance metrics from above, we formulate the verification problem as
\begin{equation}\label{prob:relative_verification}
  \begin{array}{lll}
  \mbox{maximize} & \mbox{(performance metric)} & \phi(z^K, z^\star, x)\\
  \mbox{subject to} & \mbox{(algorithm steps)}  & z^{k+1} \in s^k(z^{k}, x) \quad k=0,\dots,K-1\\
  & \mbox{(parameter)} & x \in \mathcal{X} \\
  & \mbox{(initial point)} & z^0 \in S\\
  & \mbox{(optimality)} & z^\star \in \argmin_{z \in \Omega(x)} f(z,x).
  \end{array}
\end{equation}
We let $\delta$ denote the optimal value of problem~\eqref{prob:relative_verification}.
In order for verification problem~\eqref{prob:relative_verification} to be well-posed, we rely on the following three assumptions.
\begin{assumption}[Existence of optimal solution]\label{assumption:existence_opt}
  For all parameters $x \in \mathcal{X}$, there exists a minimizer $z^\star(x)$ for problem~\eqref{prob:main}.
\end{assumption}
\begin{assumption}[Non-empty sets]\label{assumption:non_empty}
  The parameter set $\mathcal{X}$ and initial point set $S$ are both non-empty.
\end{assumption}
\begin{assumption}[Well-posed subproblems]\label{assumption:existence}
For every parameter $x \in \mathcal{X}$ and initialization $z^0 \in S$, 
each update subproblem that defines $z^{k+1}$ is well-posed.
That is, whenever the algorithm reaches iteration $k$ with iterate $z^k$,
the corresponding update operator satisfies $s^k(z^k,x) \neq \emptyset \hspace{2mm} \text{for all } k = 0,\ldots,K-1$.
\end{assumption}
Assumption~\ref{assumption:existence_opt} is mild and amounts to requiring that the parametric problem~\eqref{prob:main} always admits at least one minimizer. 
Assumption~\ref{assumption:non_empty} is also mild, requiring that there is at least one feasible parameter and initial point.
Throughout the paper, we assume that these two assumptions hold. 

Assumption~\ref{assumption:existence} requires more care.  
In many algorithmic settings this condition is natural—\eg, in the trust-region method where the constraints are convex and every admissible initial point is feasible for all parameters $x\in\mathcal{X}$.  
However, there are also important situations in which the assumption may fail, such as when downstream convex subproblems become infeasible (\eg, the final polish step in relax-round-polish).
We remark that this assumption concerns only iterates generated by the algorithm; it does not impose feasibility of the $k$-th subproblem (corresponding to $s^k(z^k,x)$) for arbitrary $z^k$, but only for those $z^k$ that arise during a valid run of the method.  

Efficiently encoding the algorithms steps and optimality conditions in problem \eqref{prob:relative_verification} is a challenge since the steps are \emph{implicitly} defined by optimization subproblems.
Moreover, the presence of tie-breakers introduces pessimism: the verification framework must account for the worst-case outcome of any tie when representing the algorithm's behavior.
In \Subsec~\ref{subsec:algorithm_steps} we show how to encode algorithm steps via their optimality conditions.
To verify suboptimality, we show how to handle the optimality condition in problem~\eqref{prob:relative_verification} in \Subsec~\ref{subsec:nonconvex_optimality}.

\subsection{Encoding algorithm steps}\label{subsec:algorithm_steps}
In this subsection, we show how to encode the algorithm steps in SCP methods as constraints in the verification problem~\eqref{prob:relative_verification}.
These building blocks are the components that we use to construct the verification problem for various algorithms.
We show how to encode convex QPs in \Subsec~\ref{subsec:convex_qp} and projections onto some non-convex sets in \Subsec~\ref{subsec:nonconvex_proj}.

\subsubsection{Convex quadratic problems}\label{subsec:convex_qp}
The main SCP steps we consider are convex QPs that take the primal and dual formulations of the form
\begin{equation}\label{prob:convex_qp}
\begin{array}{llcll}
\text{minimize}   & (1/2) u^T P u + c^T u 
& \qquad & \text{maximize} & - (1/2) u^T P u - b^T y \\[0.3em]
\text{subject to} & A u + s = b & & \text{subject to} & P u + A^T y + c = 0 \;  \\
                  & s \geq 0 & &  & y \geq 0,
\end{array}
\end{equation}
where the primal variables are $u \in \reals^{n_u}$ and $s \in \reals^{n_s}$ and the dual variable is $y \in \reals^{n_y}$.
The Karush-Kuhn-Tucker (KKT) conditions of problem~\eqref{prob:convex_qp} are~\citep{scs_quadratic}
\begin{equation*}
    A u + s = b, \quad s \geq 0, \quad P u + A^T y + c = 0, \quad y \geq 0, \quad s^T y = 0.
\end{equation*}
Since problem~\eqref{prob:convex_qp} is convex, the KKT conditions are always sufficient for optimality, and are necessary if Slater's condition holds~\citep[\Sec~5.5.3]{cvxbook}.
Since we assume that each QP is feasible (under Assumption~\ref{assumption:existence}) and the feasible region is polyhedral, Slater's condition is always satisfied~\citep[\Sec~5.2.3]{cvxbook}.
Hence, the KKT conditions are necessary and sufficient for optimality.

As is standard in the bilevel programming literature~\citep{bilevel_opt}, we avoid writing nested optimization problems in the verification problem~\eqref{prob:relative_verification} explicitly by reformulating the QPs through their KKT conditions.
In the most general case of the verification problem, the iterates $\{z^k\}_{k=0}^{K-1}$, problem parameter $x$, and the problem data of each QP $(P, A, c, b)$ (which depends on the parameter $x$ and the iterate $z^k$) all serve as decision variables.
The resulting KKT constraints introduce non-convexities (\eg, the complementarity constraint $s^T y = 0$ is a non-convex bilinear constraint).

\subsubsection{Projections onto non-convex sets}\label{subsec:nonconvex_proj}
We now consider the case where the algorithm step is a projection onto a non-convex set $\mathcal{Z} \subseteq \reals^n$ for rounding steps outlined in \Subsec~\ref{subsec:relax_round_polish}.
We denote this projection with the operator $\Pi$.

\myparagraph{Binary variables}
To formalize how the projection step onto the set $\{0,1\}^n$ can be enforced within our framework, we state the following proposition.

\begin{proposition}[Encoding binary variables]\label{prop:binary}
Let $\mathcal{Z} = \{0,1\}^n$ and suppose $u \in [0,1]^n$. 
Then
\begin{equation*}
  v = \Pi(u) \iff v \in \{0,1\}^n, \quad v - u \leq 1 / 2, \quad u - v \leq 1 / 2.
\end{equation*}
\end{proposition}
See Appendix \Subsec~\ref{proof:binary} for the proof.

\myparagraph{Sparsity constraints}
To formalize how the projection step onto the set $\{z \in \reals^n \mid \|z\|_0 \leq k\}$ can be enforced within our framework, we state the following proposition.  

\begin{proposition}[Encoding sparsity constraints]\label{prop:sparsity}
Let $\mathcal{Z} = \{z \in \reals^n \mid \|z\|_0 \leq k\}$ and suppose $u \in \reals^n$ with $w = |u|$.
Then the projection step $v = \Pi(u)$, which keeps the $k$ entries of largest magnitude and sets the rest to zero, can be equivalently encoded as
\begin{align}
    \exists \, \alpha \in \{0,1\}^n, \; \tau \geq 0 \;\; \text{s.t.} \;\;
    &\ones^T \alpha = k, \hspace{2.5mm} \alpha_i = 0 \implies w_i \leq \tau \hspace{1.5mm} \forall i, \hspace{2.5mm} \alpha_i = 1 \implies w_i \geq \tau \hspace{1.5mm} \forall i, \nonumber \\
    &\alpha_i = 0 \implies v_i = 0 \hspace{1.5mm} \forall i, \quad \alpha_i = 1 \implies v_i = u_i, \hspace{1.5mm} \forall i. \nonumber
\end{align}
\end{proposition}

See Appendix \Subsec~\ref{proof:sparsity} for the proof.
To model the logical implication constraints, we use the standard big-M technique~\citep{vielma2015mixed}, which is often supported in modern solvers~\citep{gurobi}.
We assume access to the absolute value of the vector $u$ since that variable is a natural outcome of the penalized version of the previous problem (in the relax-round-polish method).\footnote{If the absolute value is unavailable, the sparsity constraint can still easily be encoded.}

\subsection{Verifying suboptimality}\label{subsec:nonconvex_optimality}
In this subsection, we show how to verify the suboptimality of the candidate solution in cases where feasibility of the final iterate can be guaranteed.
In this subsection, we assume that the following assumption holds.
\begin{assumption}[Feasibility of final iterate]\label{assumption:feas}
  For every problem parameter $x \in \mathcal{X}$, the SCP algorithm produces a final iterate $z^K$ that lies in the feasible set $\Omega(x)$.
\end{assumption}
There are many settings where Assumption~\ref{assumption:feas} is satisfied.
\begin{itemize}[left=5pt]
  \item {\bf Unconstrained settings:} $\Omega(x) = \reals^n$, so feasibility holds automatically.  
  \item {\bf Convex feasible regions with feasibility-preserving steps:} If $\Omega(x)$ is convex for all $x \in \mathcal{X}$, the initial point is feasible for all parameters $x \in \mathcal{X}$, and each step preserves feasibility, then the final iterate remains feasible.  
  \item {\bf Non-convex feasible regions with explicit enforcement:} In some cases of simpler non-convex feasible regions, feasibility is ensured by the algorithm itself—for example, when $\Omega(x) = \{0,1\}^n$ and the final step is a rounding operation.
  \item {\bf Verified feasibility with our constraint satisfaction framework:} In some cases, none of the above three conditions may hold. 
  However, we can use the framework from \Subsec~\ref{subsec:constraint_violation} to first verify that the final iterate is always feasible after a given number of iterations.
  If feasibility can be verified, we can then solve a second verification problem for the worst-case suboptimality.
\end{itemize}

The primary difficulty in solving problem \eqref{prob:relative_verification} for worst-case suboptimality lies in enforcing the optimality condition.
In \Subsec~\ref{subsec:convex_qp} we showed how to encode the minimizer of a convex QP into the verification problem by directly writing the KKT conditions as constraints.
However, we cannot rely on this reformulation for the optimality condition in problem~\eqref{prob:relative_verification} since problem~\eqref{prob:main} is non-convex.
This means that i) the KKT conditions only encode stationary points for non-convex problems and ii) in order for the KKT conditions to be necessary for even a stationary point, additional constraint qualifications (like the linear independent constraint qualification)~\citep[Chapter~12]{NoceWrig06} must be satisfied.

It turns out that in this case where the performance metric is the suboptimality, \emph{we do not need to encode the optimality constraint}.
Instead, we relax it to a feasible constraint and formulate the problem
\begin{equation}\label{prob:simplified_verification}
  \begin{array}{lll}
  \mbox{maximize} & \mbox{(performance metric)} & f(z^K, x) - f(z^\star, x)\\
  \mbox{subject to} & \mbox{(algorithm steps)}  & z^{k+1} \in s^k(z^{k}, x) \quad k=0,\dots,K-1\\
  & \mbox{(parameter)} & x \in \mathcal{X} \\
  & \mbox{(initial point)} & z^0 \in S\\
  & \mbox{(feasibility)} & z^\star \in \Omega(x).
  \end{array}
\end{equation}
Let $\bar{\delta}$ denote the optimal value of problem~\eqref{prob:simplified_verification}.
Observe that problem~\eqref{prob:relative_verification} with objective $\phi(z^K, z^\star, x) = f(z^K, x) - f(z^\star,x)$ and problem~\eqref{prob:simplified_verification} are the same except that the optimality constraint in~\eqref{prob:relative_verification} is replaced with a feasibility constraint in~\eqref{prob:simplified_verification}.
The following theorem shows that \emph{the relaxation from problem~\eqref{prob:relative_verification} to problem~\eqref{prob:simplified_verification} is in fact tight}.

\begin{theorem}\label{thm:opt}
  Under Assumptions~\ref{assumption:existence_opt},~\ref{assumption:non_empty},~\ref{assumption:existence} and \ref{assumption:feas}, and with the performance metric defined as in \Eqn~\eqref{eq:suboptimality}, $\phi(z^K, z^\star, x) = f(z^K, x) - f(z^\star, x)$, we obtain $\delta = \bar{\delta}$.
  That is, the optimal value of problem~\eqref{prob:relative_verification} is equal to the optimal value of problem~\eqref{prob:simplified_verification}.
\end{theorem}
See Appendix \Subsec~\ref{proof:opt} for the proof.
Intuitively, the relaxation of the optimality constraint to a feasibility constraint is tight because maximizing the suboptimality \emph{forces} decision variable $z^\star$ to minimize $f(z^\star,x)$ for parameter $x$.
The tightness of the relaxation from \Thm~\ref{thm:opt} plays a pivotal role—it is what renders the verification of suboptimality doable, as it is not obvious how one would encode the optimality condition directly in problem~\eqref{prob:relative_verification}.
The approach outlined in this subsection makes it possible to obtain global worst-case guarantees on suboptimality in cases where feasibility can be certified.

\subsection{Verifying the level of constraint violation}\label{subsec:constraint_violation}
In the previous subsection, we assumed that the feasibility of the candidate solution $z^K$ could be guaranteed; yet in some cases this assumption is not realistic.
In this subsection, we discuss how to measure the worst-case constraint violation (as measured by the square of the $\ell_2$-norm of the violations) in such cases.
Assuming that the feasible region $\Omega(x)$ is defined by inequality constraints $g_i(z,x) \leq 0$, $i=1,\dots,m$, we can verify the worst-case constraint violation of the candidate solution $z^K$ by solving the verification problem
\begin{equation}\label{prob:verification_constraints}
  \begin{array}{lll}
  \mbox{maximize} & \mbox{(performance metric)} & \sum_i (\max \{g_i(z^K,x), 0\})^2 \\
  \mbox{subject to} & \mbox{(algorithm steps)}  & z^{k+1} \in s^k(z^{k}, x) \quad k=0,\dots,K-1\\
  & \mbox{(parameter)} & x \in \mathcal{X} \\
  & \mbox{(initial point)} & z^0 \in S.
  \end{array}
\end{equation}
If the optimal value of problem~\eqref{prob:verification_constraints} is zero, then the candidate solution $z^K$ is guaranteed to be feasible for any admissible parameter $x \in \mathcal{X}$.
If the optimal value is positive, then there exists a parameter $x \in \mathcal{X}$ for which the candidate solution $z^K$ is infeasible, and the optimal value quantifies the corresponding level of constraint violation.
In problem~\eqref{prob:verification_constraints} we penalize the square of the $\ell_2$-norm of violations, but we could easily adjust it to penalize other metrics (\eg, the $\ell_\infty$ or $\ell_1$ norms) of violations.

\subsection{Verifying feasibility of subproblems with linear constraints}\label{subsec:feas}
Convex subproblems arising within the SCP procedure may, in some cases, be infeasible.
In this subsection, we show how our verification framework can certify whether feasibility of a convex QP at iteration $K-1$ is guaranteed for all admissible parameters and iterates consistent with the SCP update rules.
Suppose that the convex QP at the final iteration has constraints $A^{K-1} u \leq b^{K-1}$ where the problem data $A^{K-1}$ and $b^{K-1}$ are functions of the previous iterate $z^{K-1}$ and parameter $x$: \ie, $(A^{K-1}, b^{K-1}) = \psi^{K-1}(z^{K-1}, x)$.
Feasibility of these linear constraints can be characterized using the Farkas Lemma~\citep[Section~7.3]{schrijver1986theory}.
\begin{lemma}[Farkas Lemma]\label{lemma:farkas}
    Given $A \in \reals^{m \times n}$ and $b \in \reals^m$, exactly one of the following is true:
    \begin{itemize}
        \item There exists $u \in \reals^n$ such that $Au \leq b$.
        \item There exists $y \in \reals^m$ such that $A^T y = 0$, $y \geq 0$, and $b^T y < 0$.
    \end{itemize}
\end{lemma}

Using this characterization, we embed feasibility checking into our verification framework through the following optimization problem:
\begin{equation}\label{prob:feas_subproblem}
  \begin{array}{lll}
  \mbox{maximize} & \mbox{(performance metric)} & -(b^{K-1})^T y^K\\
  \mbox{subject to} & \mbox{(algorithm steps)}  & z^{k+1} \in s^k(z^{k}, x) \quad k=0,\dots,K-2\\
  & \mbox{(parameter)} & x \in \mathcal{X} \\
  & \mbox{(initial point)} & z^0 \in S\\
  & \mbox{(polyhedral data)} & (A^{K-1}, b^{K-1}) = \psi^{K-1}(z^{K-1}, x)\\
  & \mbox{(Farkas constraints)} & y^K \geq 0, \quad (A^{K-1})^T y^K = 0,
  \end{array}
\end{equation}
where $y^K \in \reals^m$ enters as an additional decision variable.
Let $\gamma$ be the optimal value of problem~\eqref{prob:feas_subproblem}.
\begin{theorem}\label{thm:feas_subprob}
  Suppose that Assumptions~\ref{assumption:existence_opt} and~\ref{assumption:non_empty} hold and the operators $s^0,\dots,s^{K-2}$ are well-posed (\ie, Assumption~\ref{assumption:existence} holds).
  Then, the convex QP at iteration $K-1$ is feasible for all problem parameters $x \in \mathcal{X}$ and iterates consistent with SCP update rules if and only if $\gamma \leq 0$.
\end{theorem}
See Appendix \Subsec~\ref{proof:feas_subprob} for the proof which uses the Farkas Lemma~\ref{lemma:farkas}.
The significance of Theorem~\ref{thm:feas_subprob} is that it enables us to certify whether or not a downstream QP is feasible for all admissible parameters and all iterates consistent with the SCP update rules.
This is particularly useful for verifying feasibility of the polish subproblem in the relax-round-polish method described in \Subsec~\ref{subsec:relax_round_polish}.

\subsection{Inexact solves}\label{subsec:inexact}
In practice, QPs are solved using iterative QP solvers~\citep{osqp,scs_quadratic,goulart2024clarabel}, which are typically solved up to a given tolerance.
In this subsection, we describe how to handle inexact solves to problem~\eqref{prob:convex_qp} in our verification framework.
We consider two representative models of inexactness:
\begin{itemize}[left=5pt]
    \item {\bf Distance to optimality inexactness}: A common modeling assumption is that the computed solution is within an $\epsilon$-distance of the true subproblem solution~\citep[\Sec~A.2]{NoceWrig06}. 
    In this case, we include variables for the QP's optimal solution $u^{\rm exact}$ that exactly satisfies the KKT conditions and the inexact solution $u^{\rm inexact}$ with the constraint $\|u^{\rm exact} - u^{\rm inexact}\|_\infty \leq \epsilon$.
    \item {\bf KKT inexactness}: Many QP solvers, such as SCS~\citep{scs_quadratic} and OSQP~\citep{osqp}, report accuracy in terms of primal and dual residuals rather than distance to the exact solution. 
    In this case, the iterates automatically satisfy some of the optimality conditions from \Subsec~\ref{subsec:algorithm_steps} (see for example, ~\citep[\Subsec~3.3]{osqp}), but not the primal residuals $\|A u + s - b\|_\infty \leq \epsilon$ and dual residuals $\|P u + A^T y + c\|_\infty \leq \epsilon$.
    In this case, we replace the exact optimality conditions with these constraints.
\end{itemize}
We use the infinity norm since this is the standard convergence criterion in many QP solvers~\citep{scs_quadratic,osqp,perfverifyqp}.
These two models allow our framework to reliably capture the kinds of inexactness encountered in practice, ensuring that the resulting worst-case guarantees remain meaningful even when subproblems are not solved to full precision.

\section{Numerical experiments}\label{sec:numerical_experiments}
In this section, we showcase the utility of our framework to verify the performance of SCP methods with numerical examples.
We focus on the trust-region method in \Subsec~\ref{subsec:numerical_trust}, the penalized CCP in \Subsec~\ref{subsec:numerical_ccp}, the prox-linear method in \Subsec~\ref{subsec:numerical_prox_linear}, and the relax-round-polish algorithm in \Subsec~\ref{subsec:numerical_relax_round_polish}.
The code to reproduce our results is available at
\begin{equation*}
\text{\url{https://github.com/rajivsambharya/verify\_scp}}.
\end{equation*}
We solve the verification problems (which in the most general case are mixed-integer quadratically-constrained quadratic programs) to a $2 \%$ optimality gap (unless otherwise indicated) using Gurobi 12.0~\citep{gurobi}.
In some examples, we use scalability enhancements described in Appendix \Subsec~\ref{subsec:scalability} based on optimization-based bound tightening and solving the verification problems \emph{sequentially}.
We include all of the timing results to solve the verification problems in Appendix \Sec~\ref{subsec:timing}.
Our results exhibit a wide range of algorithmic behaviors—from provable convergence to a globally optimal point, to scenarios where significant progress is made but plateaus short of optimality, and cases where the algorithm immediately stalls and fails to progress.
In all examples, Assumptions~\ref{assumption:existence_opt},~\ref{assumption:non_empty}, and~\ref{assumption:existence} which ensure the existence of an optimal solution for every admissible parameter, the parameter sets and initial point sets are non-empty, and the SCP update rules are well-posed, are satisfied.

\myparagraph{Baselines}
To the best of our knowledge, no general baseline methods exist to generate worst-case guarantees for these SCP algorithms.
For validation, in all experiments we report the \emph{sample maximum}, obtained by uniformly sampling $500$ parameter instances from the set $\mathcal{X}$ (unless stated otherwise) and taking the worst-case value. 
We formulate the convex QPs in CVXPY~\citep{diamond2016cvxpy,agrawal2018rewriting} and solve them using OSQP~\citep{osqp}.

\subsection{Trust-region method}\label{subsec:numerical_trust}
In this subsection, we apply our verification framework to the trust-region method from \Subsec~\ref{subsec:trust_scp}, on box-constrained quadratic minimization in \Subsec~\ref{subsec:box_qp} and network utility maximization in \Subsec~\ref{subsec:network_utility}. 

\subsubsection{Box-constrained quadratic minimization}\label{subsec:box_qp}
We first consider the non-convex box-constrained quadratic optimization problem
\begin{equation}
  \begin{array}{ll}
  \label{prob:nonconvex_box_qp}
  \mbox{minimize} & (1 / 2) z^T P z + x^T z\\
  \mbox{subject to} & -1 \leq z \leq 1,
  \end{array}
\end{equation}
where $z \in \reals^n$ is the decision variable, $x \in \reals^n$ is the problem parameter, and the matrix $P \in \symm^{n}$ which has both positive and negative eigenvalues is shared across all problem instances.

\myparagraph{Numerical example}
We take $P = \bar{P} + \bar{P}^T$, where the entries of $\bar{P}$ are generated by sampling from a standard Gaussian distribution in an i.i.d. fashion.
In this example, we compare against two different sets for the parameters: $\mathcal{X}_1=[2,4]^n$ and $\mathcal{X}_2=[5,8]^n$.
We fix the trust-region size to be $\rho=0.2$ and solve the verification problem with a cold start $S = \{{\bf 0}_n\}$ and with a warm start $S = \{s^{\rm ws}\}$.
To calculate the warm start point $s^{\rm ws}$, we solve the non-convex problem~\eqref{prob:nonconvex_box_qp} with $x$ fixed to the point at the center of the parameter set (\eg, for $\mathcal{X}_2$, we fix $x=(6.5) \mathbf{1}_n$).
Assumption~\ref{assumption:feas} is satisfied since the initial point is always feasible and the SCP steps preserve feasibility.
We take $n=10$ in this example.
We use only $10$ samples to compute the sample maximum, as using more samples causes the curves to overlap and become difficult to distinguish.

\begin{figure}[!h]
  \centering
  \includegraphics[width=0.92\linewidth]{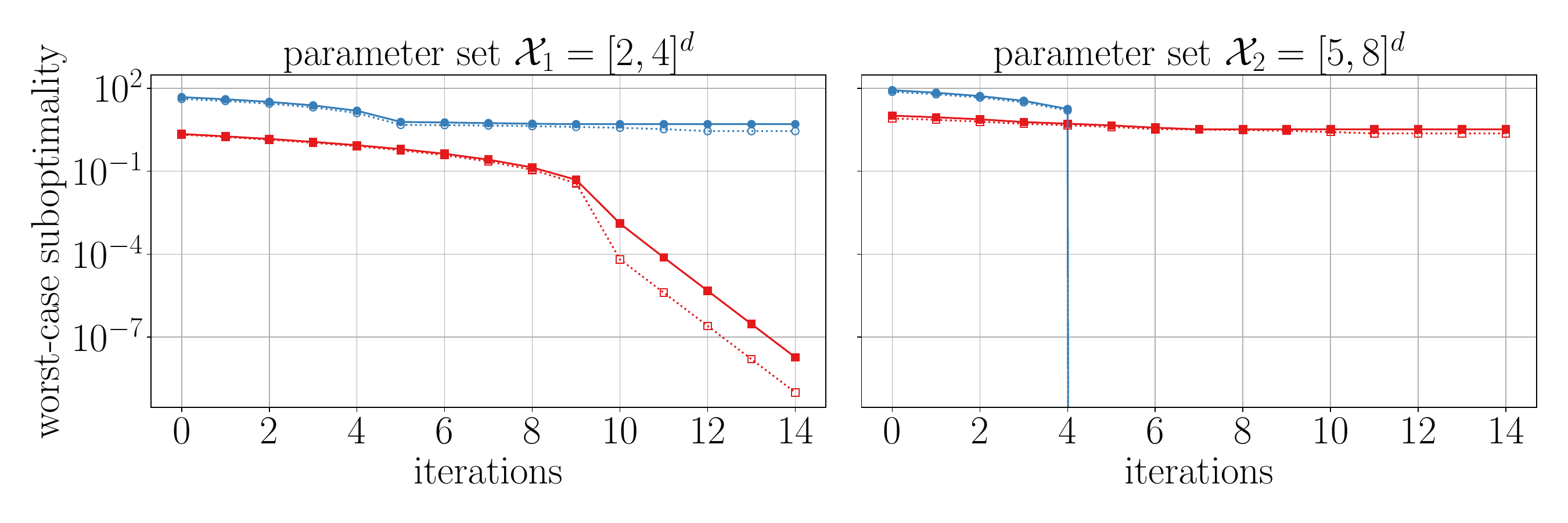}
  \\ \vspace{-3mm}
    \linesquare{228}{26}{28} verified performance: warm-started \hspace{3mm}  \linesquarehollow{228}{26}{28} sample maximum: warm-started \\
    \linecircle{55}{126}{184} verified performance: cold-started \hspace{3mm} \linecirclehollow{55}{126}{184} sample maximum: cold-started \\
    \caption{Box QP results.
    {\bf Left}: Results for $\mathcal{X}_1 = [2,4]^d$.
    The warm-started initialization is guaranteed to achieve global optimality after $14$ iterations (within tolerance $10^{-7}$), but the cold-started initialization is not guaranteed to reach a globally optimal solution.
    {\bf Right}: Results for $\mathcal{X}_2 = [5,8]^d$.
    The cold-started initialization is guaranteed to exactly achieve global optimality after $5$ iterations (the blue curve is not visible because the suboptimality is below $10^{-15}$), but the warm-started initialization is not guaranteed to reach a globally optimal solution.
    {\bf Takeaway message}: It is not obvious a priori whether warm-starting yields better results than cold-starting, or whether the trust-region method reaches global optimality. 
    Our framework provides definitive guarantees that can answer these questions by explicitly accounting for the initialization, parameter set, and the trust-region size.
    }
    \label{fig:scp_box_results}
\end{figure}

\myparagraph{Results}
We illustrate our results in Figure~\ref{fig:scp_box_results}.
In parameter set $\mathcal{X}_1$, the cold-started initialization is certifiably optimal, while the warm-started initialization is not; in parameter set $\mathcal{X}_2$, the reverse holds. 
This shows that guarantees are highly dependent on both the initialization and the parameter set.

\subsubsection{Network utility maximization}\label{subsec:network_utility}
We now turn to the resource allocation setting, modeled as a network utility maximization problem~\citep{MattingleyBoyd2010}.
We model a communication network comprising $d$ edges and $n$ flows.
Each edge $i$ has capacity $d_i$, and each flow $j$ carries a non-negative rate $z_j$.
The routing structure is encoded in a binary matrix $R \in \{0,1\}^{d \times n}$, where $R_{ij} = 1$ indicates that flow $j$ passes through edge $i$.
The total load on edge $i$ is then given by the sum of the rates of all flows that traverse it, compactly expressed as $Rz$ which must satisfy the parametric capacity constraints $Rz \leq x$.
We consider a variant with quadratic objective~\citep{fazel2005network,GuisewitePardalos1990}
\begin{equation}
  \begin{array}{ll}
  \label{prob:num}
  \mbox{maximize} & (1/2) z^T D z + d^T z\\
  \mbox{subject to} & R z \leq x, \quad 0 \leq z,
  \end{array}
\end{equation}
where $z \in \reals^n$ is the decision variable, $x \in \reals^d$ is the problem parameter, and $R \in \{0,1\}^{d \times n}$, $D \in \symm^n_+$, and $d \in \reals^n$ are problem data fixed for all instances.
Problem~\eqref{prob:num} is non-convex since the objective function being maximized is convex.

\myparagraph{Numerical example}
We set the cost problem data to be $d = {\bf 0}_{n}$ and $D = I_n$.
We generate the routing matrix $R$ as in~\citep{perfverifyqp}, where we select $1$ with $0$ for each entry i.i.d. with probability $0.5$.
We let $\mathcal{X} = [7,8]^d$.
We use a cold start (initialized from a random point in $[0,1]^d$)\footnote{We cold start from a random point instead of the zero vector because the trust-region method gets stuck at the zero vector}.
This cold start is feasible for all problems with parameter $x \in \mathcal{X}$ which means that Assumption~\ref{assumption:feas} is met.
We take $d=10$ and $n=5$ as in~\citep{perfverifyqp}.

\begin{figure}[!h]
  \centering
  \includegraphics[width=0.46\linewidth]{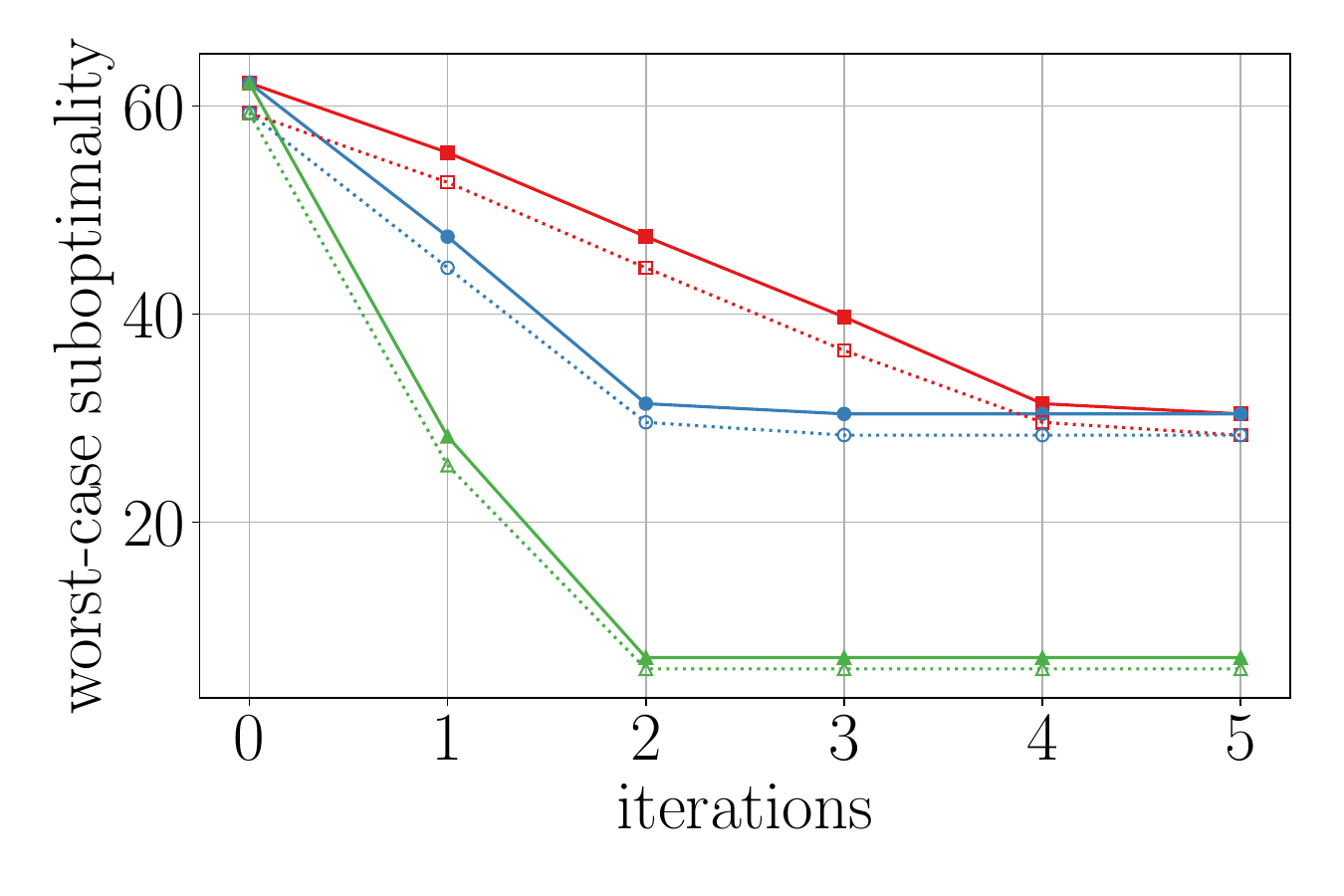}
  \\ \vspace{-3mm}
  verified performance: \linesquare{228}{26}{28} $\rho = 1$ \hspace{0.5mm} \linecircle{55}{126}{184}$\rho = 2$ \hspace{0.5mm} \lineuptri{77}{175}{74} $\rho = 5$\\
  sample maximum: \linesquarehollow{228}{26}{28} $\rho = 1$ \hspace{0.5mm} \linecirclehollow{55}{126}{184} $\rho = 2$ \hspace{0.5mm} \lineuptrihollow{77}{175}{74} $\rho = 5$
    \caption{Network utility results.
    A larger trust-region size enables the trust-region method to escape poorer local minima that the variants with a smaller trust-region size converge to—an interesting outcome, since larger steps are often associated with instability rather than improved convergence.
    }
    \label{fig:network_utility_results}
\end{figure}

\myparagraph{Results}
The results are illustrated in Figure~\ref{fig:network_utility_results}, where we test different trust-region size values $\rho$.
Our framework reveals that among the tested sizes $\{1, 2, 5\}$, the largest trust-region size $5$ provides the strongest worst-case guarantees. 
For all three cases, the worst-case guarantee plateaus within $5$ iterations.

\subsection{Penalized convex-concave procedure}\label{subsec:numerical_ccp}
In this subsection, we apply our verification framework to the penalized CCP from \Subsec~\ref{subsec:pen_convex_concave}, focusing on power converter control in \Subsec~\ref{subsec:power_converter} and the knapsack problem in \Subsec~\ref{subsec:knapsack}.

\subsubsection{Power converter control}\label{subsec:power_converter}
We consider the following optimal control problem for power electronic converters, following a similar setup from~\citep{takapoui2020simple}:
\begin{equation}
  \begin{array}{ll}
  \label{prob:optimal_control}
  \mbox{minimize} & (1/2) \sum_{t=1}^{T} (s_t - s_t^{\rm ref})^T Q_t (s_t - s_t^{\rm ref}) + (1/2) \sum_{t=1}^{T-1} u_t^T R_t u_t,\\
  \mbox{subject to} & s_{t+1} = A s_t + B u_t, \quad u_t \in \{-1,1\}^{n_u}, \quad t=0,\dots,T-1 \\
  & s_0 = s^{\rm init},
  \end{array}
\end{equation} 
where the decision variables are the states $\{s_t\}_{t=0}^{T-1}$ where $s_t \in \reals^{n_s}$ and the controls $\{u_t\}_{t=0}^{T-1}$ where $u_t \in \reals^{n_u}$.
The problem parameter is the initial state $x = s^{\rm init}$.
The dynamics matrices $A \in \reals^{n_s \times n_s}$, $B \in \reals^{n_s \times n_u}$, the reference trajectory $\{s_t^{\rm ref}\}_{t=1}^T$, and the cost matrices $Q_t \in \symm_+^{n_s}$ and $R_t \in \symm_+^{n_u}$ are fixed for all instances.
The control inputs represent the switching states of the converter, taking values in $\{-1,1\}$. 
After eliminating the state variables, problem~\eqref{prob:optimal_control} can be written in the condensed form
\begin{equation}
  \begin{array}{ll}
  \label{prob:optimal_control_condensed}
  \mbox{minimize} & (1/2) u^T P u + u^T (K x + c_0),\\
  \mbox{subject to} & u \in \{-1,1\}^{T n_u},
  \end{array}
\end{equation} 
where the decision variable $u$ is the stacked controls over the time indices.
The matrices $P \in \symm^{T n_u}$ and $K \in \reals^{T n_u \times n_s}$ and the vector $c_0 \in \reals^{T n_u}$ are derived from the problem data.
At the end of the penalized CCP, we round the solution to the nearest vector in the discrete set $\{-1,1\}^{T n_u}$ to obtain a feasible point, thus satisfying Assumption~\ref{assumption:feas} (see \Subsec~\ref{subsec:nonconvex_proj} for how to encode this step as constraints).

\myparagraph{Numerical example}
We generate the dynamics matrices $A$ and $B$ by discretizing the continuous dynamics as described in~\citep[\Subsec~3.3]{takapoui2020simple}. 
We have $n_s=4$ and $n_u=1$.
We set $Q_t$ to be the all zeroes matrix except for the bottom right entry which is $1$.
We set $R_t = 0$ and $s_t^{\rm ref} = (0, 0, 0, 1)$ for $t=0,\dots,T$.
We discretize with $\Delta t = 10^{-6}$ seconds and take $T=10$ time steps.
We set $\mathcal{X} = [-2,2]^{n_s}$ and use the hyperparameter $\tau_k = 0.2$ for the penalties in problem~\eqref{prob:pen_ccp}.
We compare the two different inexactness criteria for terminating the penalized CCP from \Subsec~\ref{subsec:inexact} for different tolerances.
We use the same warm-start initialization heuristic as described in \Subsec~\ref{subsec:box_qp}.

\begin{figure}[!h]
  \centering
  \includegraphics[width=.92\linewidth]{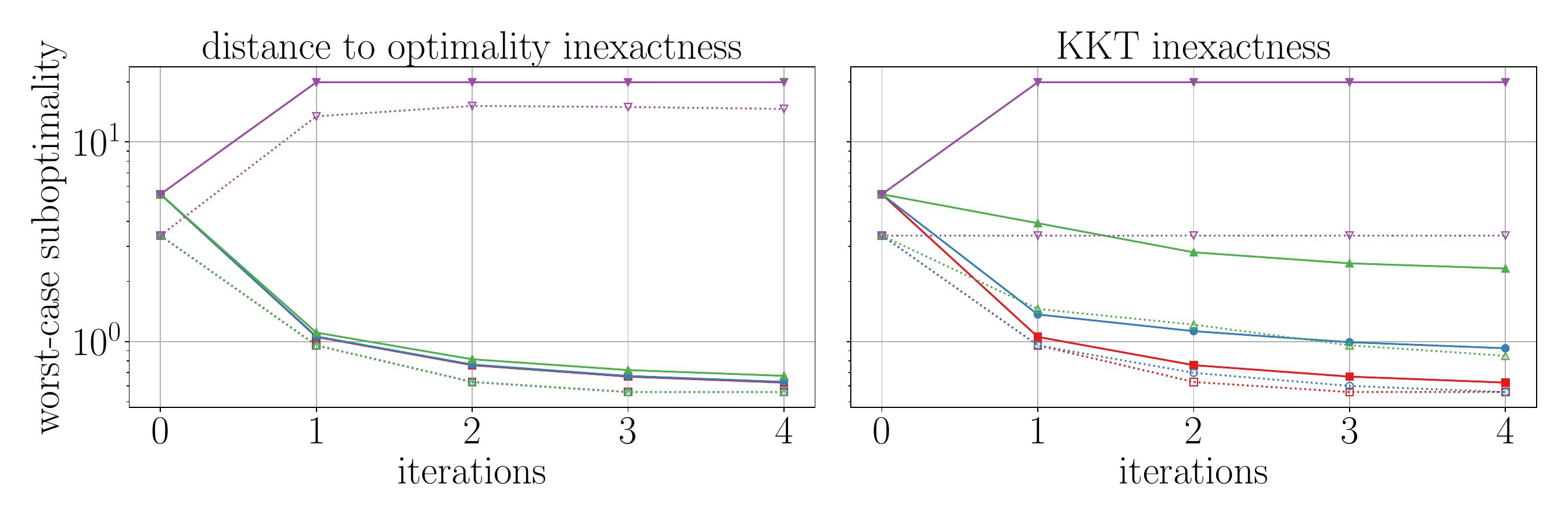}
  \\ \vspace{-3mm}
    verified performance:
    \linesquare{228}{26}{28} $\epsilon=0.0$ \hspace{0.5mm}
    \linecircle{55}{126}{184} $\epsilon=0.01$ \hspace{0.5mm}
    \lineuptri{77}{175}{74} $\epsilon=0.1$ \hspace{.5mm}\linedowntri{170}{24}{169} $\epsilon=10.0$\\ 
    sample maximum: \linesquarehollow{228}{26}{28} $\epsilon=0.0$ \hspace{0.5mm}
    \linecirclehollow{55}{126}{184} $\epsilon=0.01$ \hspace{0.5mm}
    \lineuptrihollow{77}{175}{74} $\epsilon=0.1$ \hspace{0.5mm}\linedowntrihollow{170}{24}{169} $\epsilon=10.0$
    \caption{Power converter control results.
    Left: Results for distance to optimality inexactness.
    Right: Results for KKT inexactness.
    Our framework allows us to quantitatively compare these two common stopping criteria.
    Provided a specific $\epsilon$ value, the distance to optimality inexactness criterion yields the stronger guarantees of the two. 
    }
    \label{fig:power_converter_results}
\end{figure}

\myparagraph{Results}
We illustrate our guarantees in Figure~\ref{fig:power_converter_results}.
We observe a sizeable gap between the verified worst-case performance and the sample maximum for the KKT inexactness criteria, indicating that the verification framework finds particularly adversarial inexact solutions that do not naturally arise from the OSQP solver.

\subsubsection{Knapsack}\label{subsec:knapsack}
We now consider the knapsack problem formulated as
\begin{equation}
  \begin{array}{ll}
  \label{prob:knapsack}
  \mbox{maximize} & x^T z\\
  \mbox{subject to} & a^T z \leq b, \quad z \in \{0,1\}^{n},
  \end{array}
\end{equation} 
where $z \in \reals^n$ is the decision variable, $x \in \reals^n$ is the problem parameter, and $a \in \reals_{++}^n$ and $b \in \reals_{++}$ are problem data that are fixed for all problem instances.
The goal is to select a subset of items (given by binary variables $z$) that maximizes the total value $x^T z$ while ensuring the total occupied space $a^T z$ does not exceed capacity $b$.

\myparagraph{Numerical example}
We set $n=10$, randomly sample the vector $a$ from a uniform distribution in $[0,1]$ in an i.i.d. fashion, and set $b$ to be half the sum of the entries of $a$.
We set the penalty hyperparameters with $\tau_k = \tau_0 \kappa^k$ and test across $\tau_0 \in \{0.01, 1, 100\}$ and $\kappa \in \{1, 2\}$.
We initialize from the vector $(0.5)\ones_n$ (the middle of the set $\{0,1\}^n$) and take the parameter set $\mathcal{X} = [5,7]^n$.

\begin{figure}[!h]
  \centering
  \includegraphics[width=.92\linewidth]{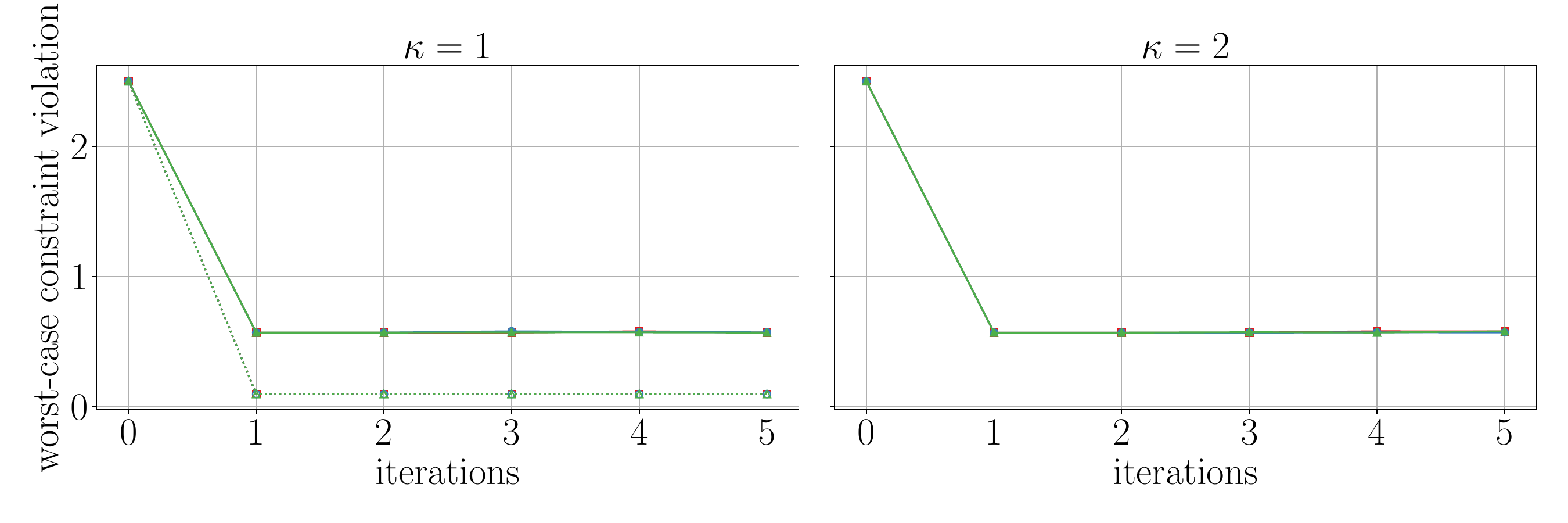}
  \\ \vspace{-3mm}
    verified performance:
    \linesquare{228}{26}{28} $\tau_0=0.01$ \hspace{0.5mm}
    \linecircle{55}{126}{184} $\tau_0=1$ \hspace{0.5mm}
    \lineuptri{77}{175}{74} $\tau_0=100$ \hspace{.5mm}\\ 
    sample maximum: \linesquarehollow{228}{26}{28} $\tau_0=0.01$ \hspace{0.5mm}
    \linecirclehollow{55}{126}{184} $\tau_0=1$ \hspace{0.5mm}
    \lineuptrihollow{77}{175}{74} $\tau_0=100$ \hspace{0.5mm}
    \caption{Knapsack results.
    All of the verification curves overlap with each other across all $\tau_0$ and $\kappa$ values.
    In all cases, the algorithm makes progress for $1$ iteration and then stalls.
    When $\kappa=2$, the final penalty value becomes very large (which is designed to prioritize feasibility~\citep[\Subsec~3.1]{lipp2016variations}), yet the algorithm still fails to find a feasible solution.
    This example shows how our framework can reveal cases in which the SCP algorithm is ineffective.
    Finally, for all $\kappa=2$ settings, the verified performance matches the sample maximum exactly.
    }
    \label{fig:knapsack_results}
\end{figure}

\myparagraph{Results}
We showcase our bounds on the level of constraint violation in Figure~\ref{fig:knapsack_results}.
Across all choices of hyperparameters $\tau_0$ and $\kappa$, the penalized CCP fails to make any progress beyond the first iteration, and no feasible solution can be verified in any case.

\subsection{Prox-linear method}\label{subsec:numerical_prox_linear}
In this subsection, we apply our verification framework to the prox-linear method from \Subsec~\ref{subsec:prox_linear} for phase retrieval in \Subsec~\ref{subsec:phase_retrieval}. 

\subsubsection{Phase retrieval}\label{subsec:phase_retrieval}
We consider the phase retrieval problem
\begin{equation}
  \begin{array}{ll}
  \label{prob:phase_retrieval}
  \mbox{minimize} & (1 / d) \sum_{i=1}^d |(a_i^T z)^2 - x_i|,
  \end{array}
\end{equation}
where $z \in \reals^n$ is the decision variable, $x \in \reals^d$ is the problem parameter, and the dictionary vectors $\{a_i\}_{i=1}^d$ (where $a_i \in \reals^n$) are fixed across all instances.
Phase retrieval seeks to reconstruct a signal $z$ from quadratic measurements, which arises in applications such as optics, imaging, and crystallography~\citep{drusvyatskiy2017proximal}.

\myparagraph{Numerical example}
We consider a small phase retrieval problem where $d=15$ and $n=5$.
We use hyperparameter $\rho=1$ and compare against two different parameter sets: $\mathcal{X}_1 = [6,7]^d$ and $\mathcal{X}_2 = [6.5,7]^d$.
We generate the entries of the $a_i$ vectors i.i.d. from a standard Gaussian distribution with probability $0.25$ and otherwise set it to $0$.
We use the same warm-start initialization heuristic as described in \Subsec~\ref{subsec:box_qp}.
Since problem~\eqref{prob:phase_retrieval} is unconstrained, Assumption~\ref{assumption:feas} is met.

\begin{figure}[!h]
  \centering
  \includegraphics[width=0.46\linewidth]{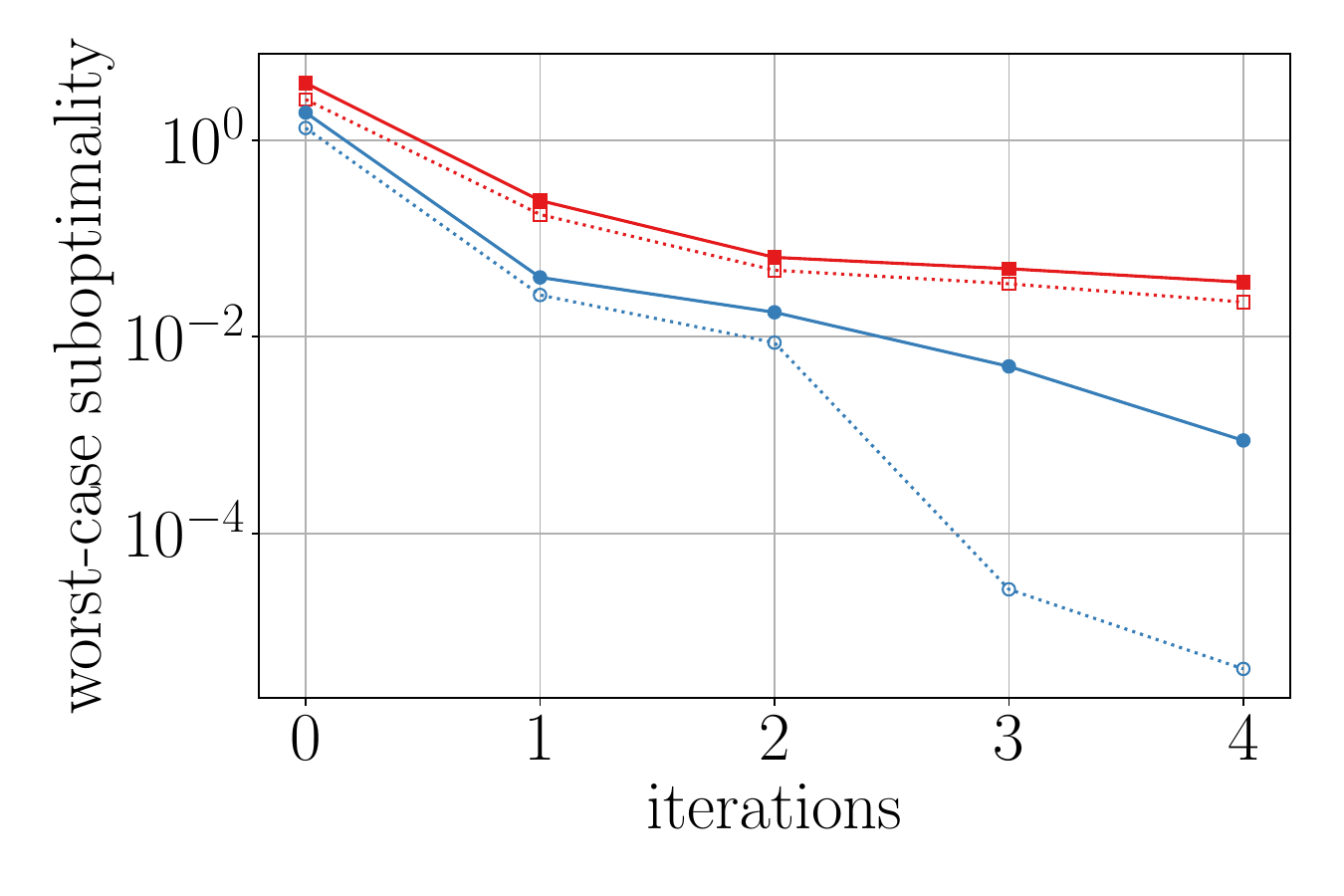}
  \\ \vspace{-3mm}
  \linesquare{228}{26}{28} verified performance: set $\mathcal{X}_1$ \hspace{3mm}  \linesquarehollow{228}{26}{28} sample maximum: set $\mathcal{X}_1$ \\
    \linecircle{55}{126}{184} verified performance: set $\mathcal{X}_2$ \hspace{3mm} \linecirclehollow{55}{126}{184} sample maximum: set $\mathcal{X}_2$
    \caption{Phase retrieval results.
    With the smaller parameter set $\mathcal{X}_2$, we can certify global optimality (up to tolerance $0.001$) after $4$ iterations.
    The worst-case guarantee for the larger parameter set $\mathcal{X}_1$ is significantly worse. 
    For $K=4$ and $\mathcal{X}_2$, Gurobi fails to reach a $2\%$ optimality gap within the $2$ hour time limit, so we report the best upper bound obtained.
    }
    \label{fig:phase_retrieval_results}
\end{figure}

\myparagraph{Results}
The results for this example are illustrated in Figure~\ref{fig:phase_retrieval_results}.
For the smaller parameter set, we can obtain significantly stronger worst-case guarantees.

\subsection{Relax-round-polish}\label{subsec:numerical_relax_round_polish}
In this subsection, we apply our framework to verify the relax-round-polish method from \Subsec~\ref{subsec:relax_round_polish} on sparse coding in \Subsec~\ref{subsec:sparse_coding} and hybrid vehicle control in \Subsec~\ref{subsec:hybrid_vehicle}.

\subsubsection{Sparse coding}\label{subsec:sparse_coding}
We consider the cardinality-constrained version of a sparse coding problem
\begin{equation}
  \begin{array}{ll}
  \label{prob:lasso}
  \mbox{minimize} & (1 / 2) \|A z - x\|_2^2\\
  \mbox{subject to} & \|z\|_0 \leq k,
  \end{array}
\end{equation}
where $z \in \reals^n$ is the decision variable, $x \in \reals^d$ is the problem parameter, and $A \in \reals^{d \times n}$ is problem data that is fixed for all instances.
In the relax-round-polish method, recall that our relax step consists in solving the lasso problem $\min_z (1 / 2) \|A z - x\|_2^2 + \lambda \|z\|_1$ with hyperparameter $\lambda$.
The subsequent round step retains only the $k$ largest entries of $z$ in magnitude, setting all others to zero.
Finally, the polish step fixes this support and solves a least-squares problem over the selected variables.
We use our verification framework to upper bound the suboptimality.

\myparagraph{Numerical example}
We follow a similar setup from~\citep{ranjan2024exact} to generate a random $A$ matrix.
With probability $0.2$, we sample the entries of $A \in \reals^{d \times n}$ i.i.d. from a standard normal distribution and otherwise set each entry to $0$.
We then normalize the columns to have unit norm.
We take the set $\mathcal{X} = [5,8]^d$ and use $d=20$, $n=15$, and $k=5$.
We then solve the verification problem for a range of $\lambda$ values spaced evenly on a log-scale between $0.01$ and $100$.

\begin{figure}[!h]
  \centering
  \includegraphics[width=0.46\linewidth]{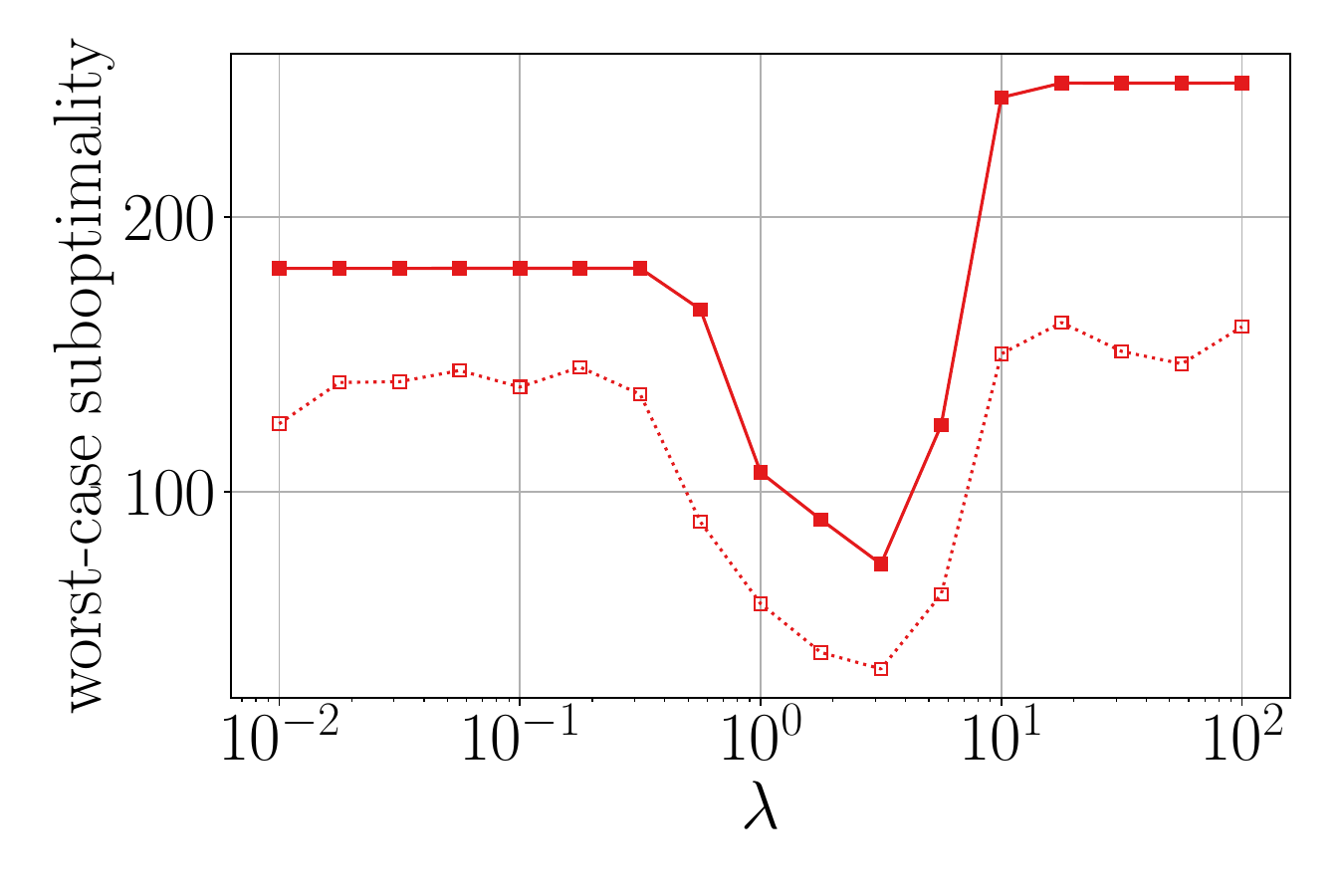}
    \\ \vspace{-3mm}
    \linesquare{228}{26}{28} verified performance \hspace{4mm}  \linesquarehollow{228}{26}{28} sample maximum 
    \vspace{-1mm}
    \caption{Sparse coding results.
    Among the set of hyperparameter options, $\lambda=3.16$ yields the lowest worst-case suboptimality. 
    We observe a sizeable gap between the verified bound and the sample maximum, which we attribute to the method's combinatorial structure.
    By selecting an adversarial parameter, the verification procedure can (i) induce ties in the relaxation that make multiple rounding outcomes equally valid, and (ii) break those ties in the most adverse way to maximize suboptimality.
    In contrast, the sampling procedure breaks such ties at random, leading to milder observed outcomes.
    }
    \label{fig:lasso_results}
\end{figure}

\myparagraph{Results}
Figure~\ref{fig:lasso_results} shows the results for different $\lambda$ values.
The verification results indicate that a globally optimal solution is not guaranteed for any value of~$\lambda$. 
Moreover, our framework can be used to select $\lambda$ to optimize worst-case performance.

\subsubsection{Hybrid vehicle control}\label{subsec:hybrid_vehicle}
We study a hybrid vehicle control problem that captures the task of meeting a power demand while managing battery energy, limiting engine switching, and enforcing operational constraints.
Formally, we consider
\begin{equation}
  \begin{array}{ll}
  \label{prob:fuel_cell}
  \mbox{minimize} & \eta (E_T - E^{\rm max})^2 + \sum_{t=0}^{T-1} (\alpha (P_t^{\rm eng})^2 + \beta P^{\rm eng}_t + \gamma z_t + c (z_{t+1} - z_t)_+)\\
  \mbox{subject to} & E_{t+1} = E_t - \tau P^{\rm batt}_t, \quad t=0,\dots,T-1 \\
  & P^{\rm batt}_t + P^{\rm eng}_t = P^{\rm load}_t, \quad t=0,\dots,T-1\\
  & E^{\rm min} \leq E_t \leq E^{\rm max}, \quad 0 \leq P_t^{\rm eng} \leq z_t P^{\rm max} \quad t=0,\dots,T-1\\
  & E_0 = E^{\rm init}, \quad z_{-1} = z^{\rm prev} \\
  & z_t \in \{0,1\}, \quad t=0,\dots,T-1,
  \end{array}
\end{equation} 
where the decision variables are $\{P^{\rm batt}_t\}_{t=0}^{T-1}$, $\{P^{\rm eng}_t\}_{t=0}^{T-1}$, $\{E_t\}_{t=0}^{T}$, and $\{z_t\}_{t=-1}^{T-1}$.
The problem parameter is $x=(E^{\rm init}, \{P^{\rm load}_t\}_{t=0}^{T-1}, z^{\rm prev}) \in \reals \times \reals^{T} \times \{0,1\}$.
The scalars $\alpha$ and $\beta$ penalize engine use, $\eta$ encourages the final energy level to be high, $\gamma$ penalizes keeping the engine on, and $c$ penalizes switching the engine on and off.
The scalars $E^{\rm min}$, $E^{\rm max}$, and $P^{\rm max}$ give the operating limits of the battery and engine.
The horizon length is given by $T$ and the sampling time by $\tau$.
In this example, we focus on certifying the feasibility of the polish step in relax-round-polish.

\myparagraph{Numerical example}
We take $\alpha=1$, $\beta=10$, $\gamma=1.5$, $\eta=2$, $c=1$, $P^{\rm max} = 1$, $E^{\rm min} = 5$, $E^{\rm max} = 10$, $\tau=2$, and $T=5$.
We take the parameter set to be $\mathcal{X} = [E^{\rm mid} - \Delta E, E^{\rm mid} + \Delta E] \times [P^{\rm load}_{\rm lb}, P^{\rm load}_{\rm ub}]^T \times\{0,1\}$, where $E^{\rm mid} = 7.5$.
We conduct three different experiments with $\Delta E \in \{0, 0.5, 1\}$.
Within each experiment, we conduct a grid search on $P^{\rm load}_{\rm lb}$ and $P^{\rm load}_{\rm ub}$ (for $P^{\rm load}_{\rm lb} \leq P^{\rm load}_{\rm ub}$).
In this example, we solve the verification problems to the optimality gap tolerance of $10^{-9}$.

\begin{figure}[!h]
  \centering
  \includegraphics[width=0.92\linewidth]{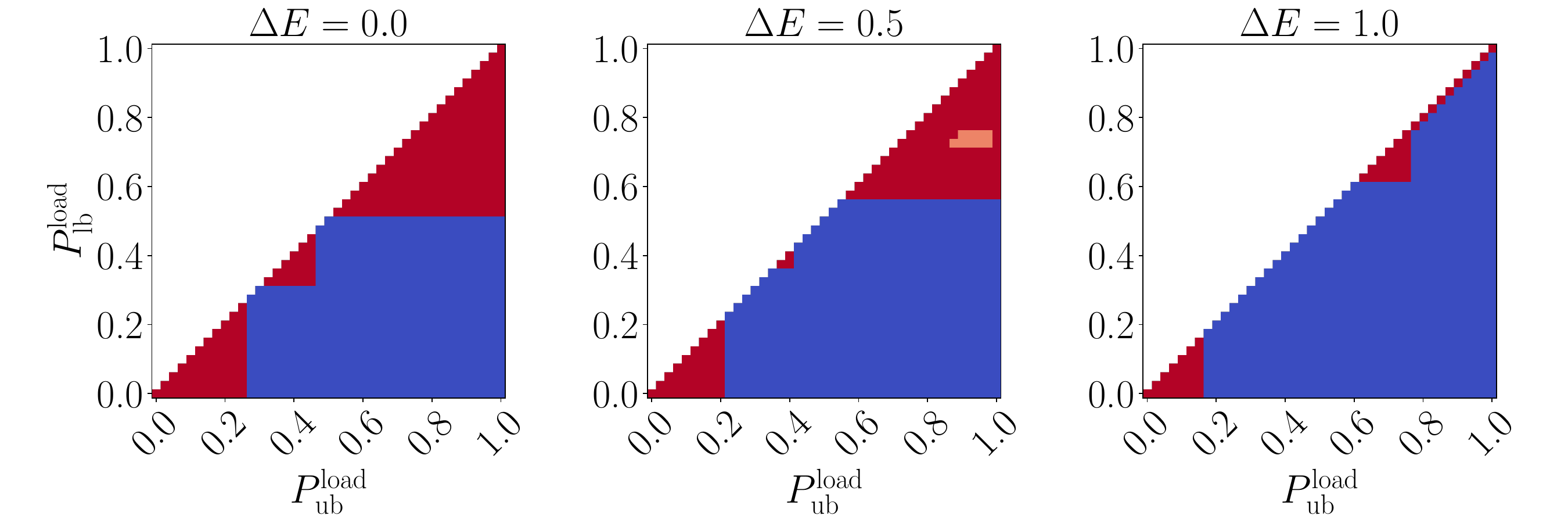}
    \\ 
    \cblock{180}{1}{1} feasibility certificate found \hspace{4mm}  \cblock{88}{86}{200} infeasibility certificate found \\
    \cblock{240}{150}{140} time limit reached, but feasibility verified from other verification problems
    \caption{Hybrid vehicle results.
    Left: $\Delta E = 0.0$.
    Middle: $\Delta E = 0.5$.
    Right: $\Delta E = 1.0$.
    As we go from i) right-to-left across the three plots, ii) right-to-left within any plot, or iii) down-to-up within any plot, the size of the parameter set $\mathcal{X}$ strictly decreases, so it is strictly easier to find the feasibility certificate. 
    We use this fact to verify feasibility in cases, where the time limit is reached.
    }
    \label{fig:hybrid_vehicle_results}
\end{figure}

\myparagraph{Results}
Figure~\ref{fig:hybrid_vehicle_results} illustrates the certificates returned for different values of $\Delta E$, $P_{\rm lb}^{\rm load}$, and $P_{\rm ub}^{\rm load}$.
In a few cases, the solver hits the time limit, but feasibility can still be inferred from neighboring certificates.
The plots show that our method cleanly separates feasible and infeasible regions for different parameter sets $\mathcal{X}$.

\section{Conclusion}\label{sec:conclusion}
We introduce a verification framework to certify the worst-case performance of sequential convex programming methods for parametric non-convex optimization.
The verification problem is formulated as an optimization problem that maximizes a performance metric over a given parameter set subject to constraints representing the SCP steps.
Our framework is general and applies to many SCP algorithms and extends naturally to inexact subproblem solves. 
Applications in control, signal processing, and operations research demonstrate that the framework provides worst-case guarantees in settings where they are otherwise unavailable.

There are several interesting directions for future work.
First, it would be interesting to extend our framework to handle semidefinite program relaxations of non-convex optimization problems and low-rank rounding schemes (as done in~\citep{goemans1995improved}).
Second, adapting the framework to analyze randomized rounding schemes would enable certification of algorithms that incorporate stochastic elements.
Third, improving the scalability of the approach remains an important challenge, for example through more efficient numerical implementations or problem-specific simplifications.

\section*{Acknowledgements} NM is supported in part by NSF Award SLES-2331880, NSF CAREER award ECCS-2045834, and AFOSR Award FA9550-24-1-0102.

\bibliography{bibliographynourl}



\appendix

\section{Proofs}

\subsection{Proof of Proposition~\ref{prop:binary}}\label{proof:binary}
\begin{proof}
Consider a single coordinate $i \in \{1,\dots,n\}$.
Since the argument holds independently for each coordinate, the claim follows for all $i = 1,\dots,n$.

$(\implies)$ If $v_i = \Pi(u_i)$, then by definition $v_i = 0$ when $u_i \leq 1/2$ and $v_i = 1$ when $u_i \geq 1/2$, so the inequalities are satisfied.

$(\impliedby)$ If $v_i \in \{0,1\}$ and the inequalities hold, then $u_i \leq 1/2$ forces $v_i = 0$, and $u_i \geq 1/2$ forces $v_i = 1$, which is exactly $v_i = \Pi(u_i)$.



\end{proof}

\subsection{Proof of Proposition~\ref{prop:sparsity}}\label{proof:sparsity}
\begin{proof}
Let $u \in \reals^n$ and $w = |u|$. Denote by $|w|_{[1]} \ge \cdots \ge |w|_{[n]}$ the order statistics of $\{w_i\}_{i=1}^n$.

\myparagraph{($\implies$) From $v=\Pi(u)$ to the encoding}
Let $\tau = w_{[k]}$ denote the $k$-th largest value among $\{w_i\}_{i=1}^n$. We consider two cases.

\emph{Case 1: $\|u\|_0 \ge k$.}
By definition of the Euclidean projection onto $\mathcal{Z}$, $v$ keeps exactly $k$ indices of (one choice of) largest magnitude and zeros the rest. 
Let $S = \{i_{[1]},\dots,i_{[k]}\}$; \ie, $S$ contains the indices of the $k$ largest-magnitude entries of $u$.
Set $\alpha_i = \mathcal{I}_S(i)$. 
Observe that $\alpha^T \ones_n = k$.
Then $\min_{i\in S} w_i \ge \tau \ge \max_{i\notin S} w_i$ and $v_i = u_i$ if $\alpha_i=1$ and $v_i=0$ if $\alpha_i=0$, so the implications in the encoding hold.

\emph{Case 2: $\|u\|_0 < k$.}
Here $\Pi(u)=u$. Let $S$ contain all nonzero indices of $u$ (so $|S|=\|u\|_0$) and add any $(k-\|u\|_0)$ indices with $w_i=0$ to reach $|S|=k$. Set $\alpha_i=\mathcal{I}_S(i)$ and choose $t=0$. Then for $i\in S$, $w_i\ge 0 = t$ and $v_i=u_i$; for $i\notin S$, $w_i\le 0 = t$ and $v_i=0$. Hence the encoding holds.

In both cases (with ties handled by the freedom in choosing $T$ and $t$), the existence of $(\alpha,t)$ satisfying the stated conditions follows from $v=\Pi(u)$.

\myparagraph{($\impliedby$) From the encoding to $v=\Pi(u)$}
Conversely, suppose there exist $\alpha \in \{0,1\}^n$ and $\tau \ge 0$ such that
\[
\ones^\top \alpha = k,\hspace{2mm}
\alpha_i=1 \Rightarrow w_i \ge \tau,\hspace{2mm}
\alpha_i=0 \Rightarrow w_i \le \tau,\hspace{2mm}
\alpha_i=1 \Rightarrow v_i = u_i,\hspace{2mm}
\alpha_i=0 \Rightarrow v_i=0.
\]
Let $S = \{ i \mid \alpha_i=1\}$; then $|S|=k$, $v_i=u_i$ for $i\in S$, and $v_i=0$ for $i\notin S$. The threshold conditions imply $\min_{i \in S} w_i \;\ge \tau \ge \max_{i \notin S} w_i$, so every selected index has magnitude at least as large as any unselected index. 
Hence $S$ indexes (one choice of) the $k$ largest magnitudes of $u$, and $v$ keeps exactly those entries while zeroing the rest. 
This is precisely the projection $v=\Pi(u)$ onto $\mathcal{Z}$ (up to tie-breaking).
\end{proof}

\subsection{Proof of \Thm~\ref{thm:opt}}\label{proof:opt}
\begin{proof}
  For notational convenience, we denote the optimal solutions in problem~\eqref{prob:relative_verification} with a subscript ``${\rm opt}$'' (\eg, $z^\star_{\rm opt}$) and the optimal values for the decision variables in problem~\eqref{prob:simplified_verification} with a subscript ``${\rm feas}$'' (\eg, $z^\star_{\rm feas}$).
  Under Assumptions~\ref{assumption:existence_opt},~\ref{assumption:non_empty}, and~\ref{assumption:existence} both problems are feasible.
  Under Assumption~\ref{assumption:feas}, $z^K_{\rm opt}\in \Omega(x_{\rm opt})$ and $z^K_{\rm feas} \in \Omega(x_{\rm feas})$.
  It is obvious that $\delta \leq \bar{\delta}$ since the optimality constraint implies the feasibility constraint.
  We must now prove that $\delta \geq \bar{\delta}$ under Assumption~\ref{assumption:existence} and choice of $\phi$.
  To do so, it is sufficient to prove that $z^\star_{\rm feas}$ is optimal for problem~\eqref{prob:main} with parameter $x_{\rm feas}$.
  Suppose it was not true. 
  Then there exists a point $\bar{z} \in \reals^n$ that is feasible in problem~\eqref{prob:main} with parameter $x_{\rm feas}$ such that $f(\bar{z},x_{\rm feas}) \leq f(z^\star_{\rm feas},x_{\rm feas})$.
  But if this was true, then $z^\star_{\rm feas}$ cannot be optimal to~\eqref{prob:simplified_verification} since replacing it with $\bar{z}$ improves the objective.
  This finishes the proof by contradiction.
\end{proof}

\subsection{Proof of \Thm~\ref{thm:feas_subprob}}\label{proof:feas_subprob}
\begin{proof}
Under Assumptions~\ref{assumption:existence_opt},~\ref{assumption:non_empty}, and~\ref{assumption:existence}, problem~\eqref{prob:feas_subproblem} is feasible.
By the Farkas Lemma (Lemma~\ref{lemma:farkas}), the constraint set $\{u \mid A^{K-1} u \le b^{K-1}\}$ is non-empty if and only if there is no certificate $y^K$ satisfying $(A^{K-1})^\top y^K = 0, y^K \ge 0, (b^{K-1})^\top y^K < 0$.
In the verification problem~\eqref{prob:feas_subproblem}, we maximize $-(b^{K-1})^\top y^K$ over all parameters $x \in \mathcal{X}$, iterates consistent with the SCP update rules, and dual variables $y^K$ satisfying $
(A^{K-1})^\top y^K = 0, y^K \ge 0$.
Hence, the optimal value $\gamma$ is the largest attainable value of $-(b^{K-1})^\top y^K$.
If $\gamma > 0$, an infeasibility certificate exists. 
Conversely, if $\gamma \leq 0$, no infeasibility certificate can be constructed, and the QP is always feasible.
\end{proof}

\section{Scalability improvements}\label{subsec:scalability}
We now present scalability enhancements that enable solving the verification problem over a greater number of iterations.
We use this strategy for two examples: phase retrieval and sparse coding.

\myparagraph{Optimization-based bound tightening}
We take advantage of optimization-based bound tightening (OBBT)~\citep{gleixner2017three} to tighten the feasible region of the verification problem.
The key idea is to tighten the lower and upper bounds on each decision variable by solving auxiliary optimization problems that minimize or maximize that variable subject to the verification constraints.
Although this procedure can be applied iteratively to further tighten the bounds, we perform only a single pass in our examples.
Each bound-tightening subproblem is solved with a time limit of $5$ seconds.

\myparagraph{Sequential solution strategy}
For algorithms with a variable number of iterations (\ie, all methods except relax-round-polish), we solve the verification problems sequentially, starting from iteration $0$ and proceeding up to $K$.
The bounds obtained from the OBBT procedure (if used) above serve as valid upper and lower bounds for each variable and are reused across iterations to improve solver efficiency.

\section{Timing results}\label{subsec:timing}
\begin{figure}[!h]
  \captionsetup{skip=4pt}
  \centering
  \includegraphics[width=0.88\linewidth]{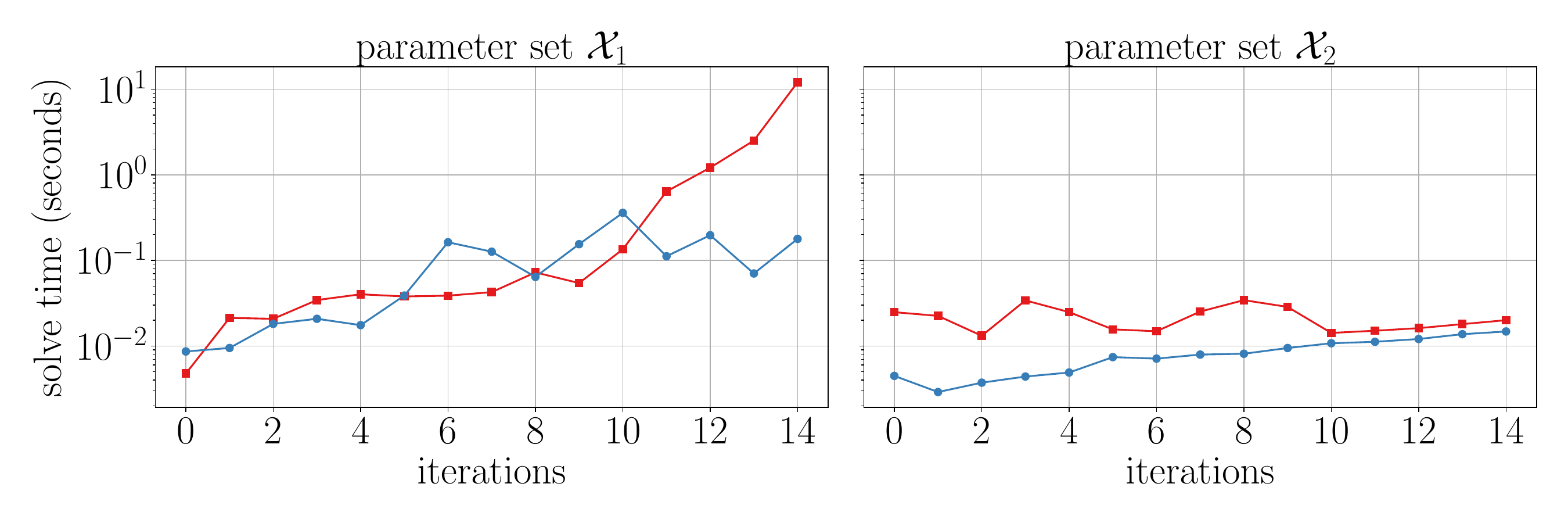}
  \\ \vspace{-3mm}
    {\small \linecircle{55}{126}{184} cold-started  \hspace{2mm}
    \linesquare{228}{26}{28} warm-started  \hspace{2mm}}\\
    \caption{\vspace{-2mm}Timing results for the box QP from \Subsec~\ref{subsec:box_qp}.
    }
    \label{fig:box_qp_timing_results}
\end{figure}

\begin{figure}[!h]
  \captionsetup{skip=4pt}
  \centering
  \includegraphics[width=0.44\linewidth]{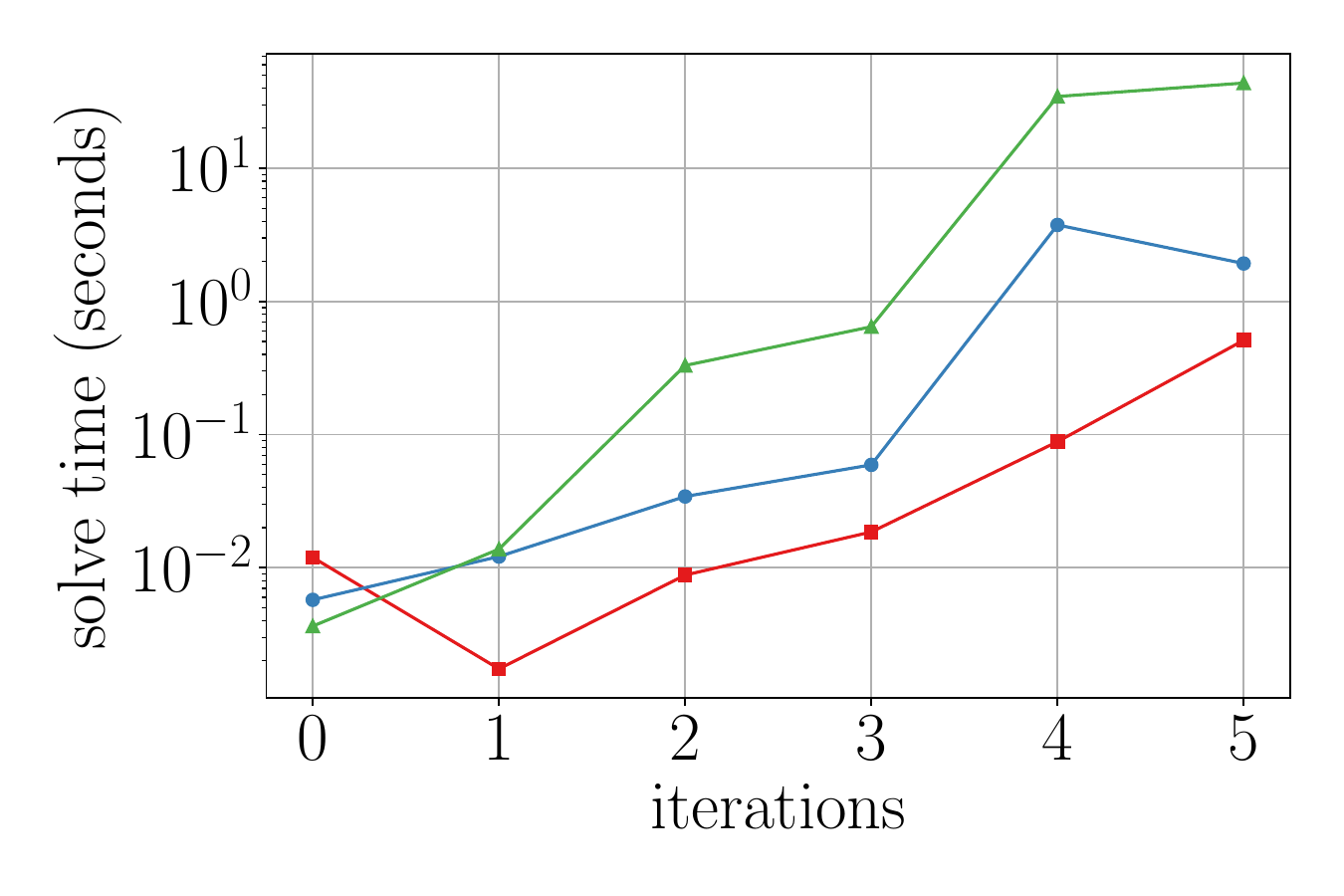}
  \\ \vspace{-3mm}
  {\small \linesquare{228}{26}{28} $\rho = 1.0$ \hspace{2mm}
    \linecircle{55}{126}{184} $\rho = 2.0$ \hspace{2mm}
    \lineuptri{77}{175}{74} $\rho = 5.0$ \hspace{2mm}}\\
    \caption{\vspace{-2mm}Timing results for the network utility maximization problem from \Subsec~\ref{subsec:network_utility}.
    }
    \label{fig:network_utility_timing_results}
\end{figure}

\begin{figure}[!h]
  \captionsetup{skip=4pt}
  \centering
  \includegraphics[width=0.88\linewidth]{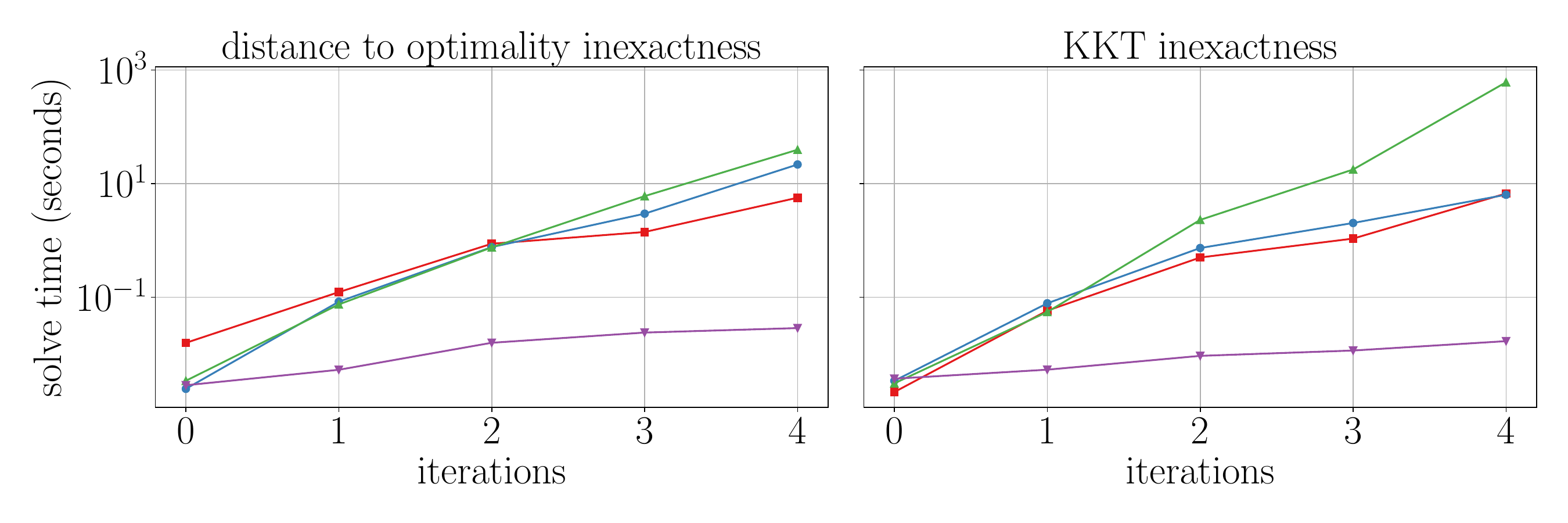}
  \\ \vspace{-3mm}
  {\small
  \linesquare{228}{26}{28} $\epsilon = 0.0$ \hspace{2mm}
    \linecircle{55}{126}{184} $\epsilon = 0.01$ \hspace{2mm}
    \lineuptri{77}{175}{74} $\epsilon = 0.1$ \hspace{2mm}
    \linedowntri{170}{24}{169} $\epsilon = 10.0$ \hspace{2mm}}\\
    \caption{\vspace{-2mm}Timing results for the power converter problem from \Subsec~\ref{subsec:power_converter}.
    }
    \label{fig:power_converter_timing_results}
\end{figure}

\begin{figure}[!h]
  \captionsetup{skip=4pt}
  \centering
  \includegraphics[width=0.88\linewidth]{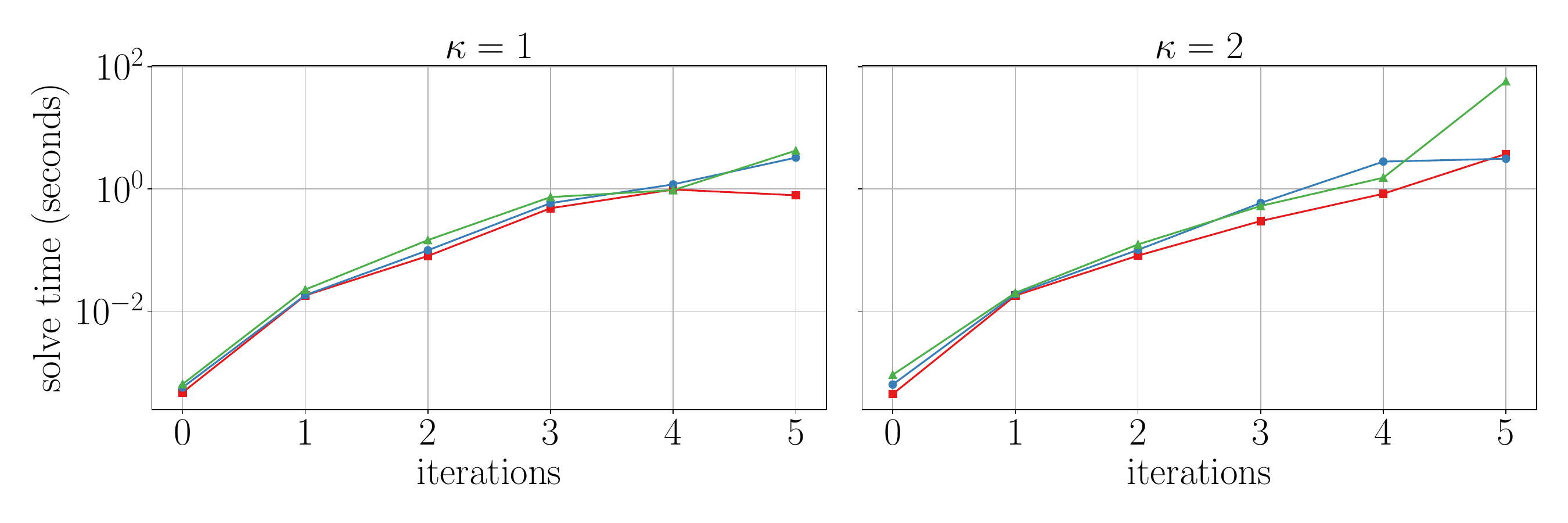}
  \\ \vspace{-3mm}
  {\small
  \linesquare{228}{26}{28} $\tau_0 = 0.01$ \hspace{2mm}
    \linecircle{55}{126}{184} $\tau_0 = 1.0$ \hspace{2mm}
    \lineuptri{77}{175}{74} $\tau_0 = 100.0$ \hspace{2mm}}\\
    \caption{\vspace{-2mm}Timing results for the knapsack problem from \Subsec~\ref{subsec:knapsack}.
    }
    \label{fig:knapsack_timing_results}
\end{figure}

\begin{figure}[!h]
  \captionsetup{skip=4pt}
  \centering
  \includegraphics[width=0.44\linewidth]{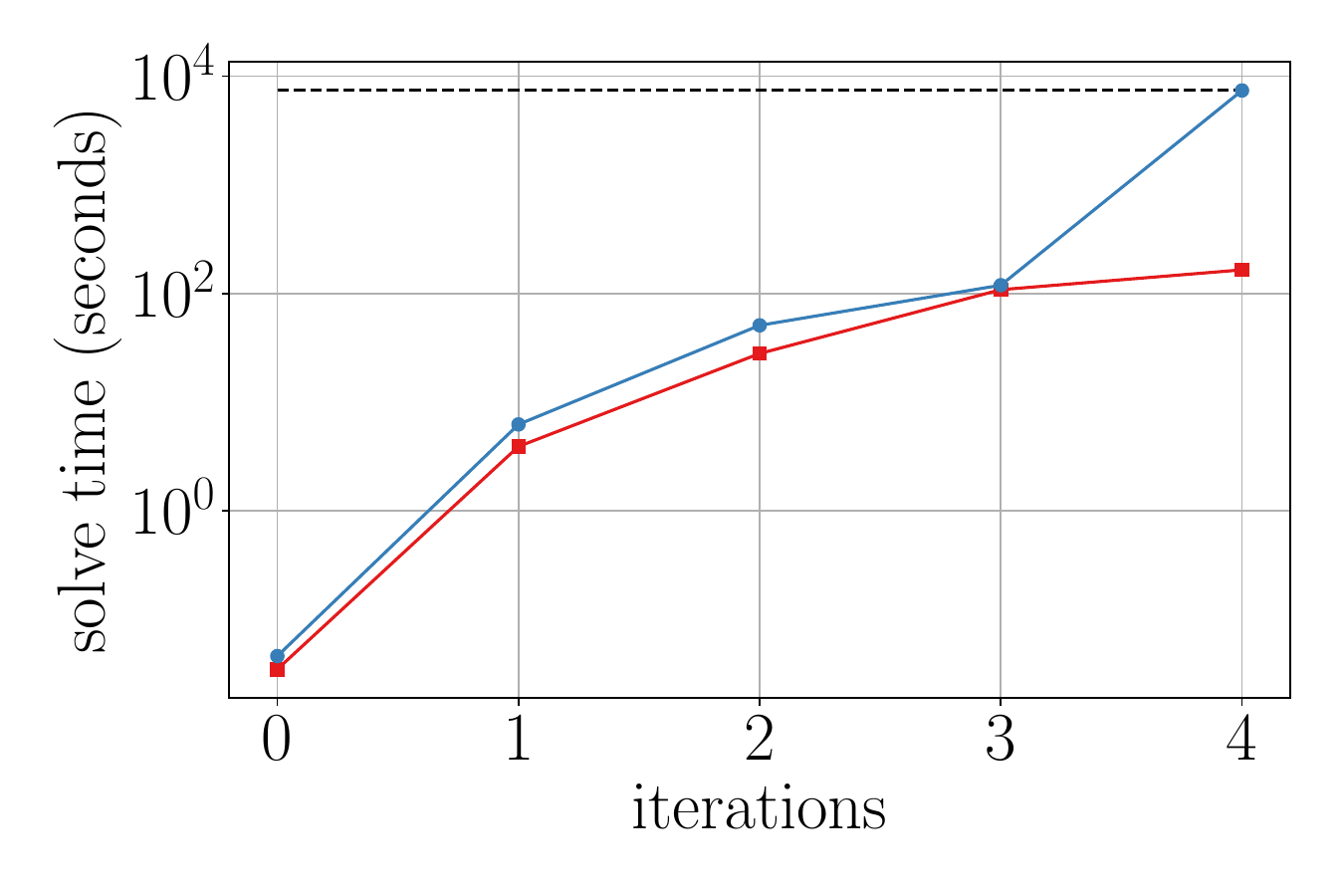}
  \\ \vspace{-2mm}
  {\small
   \linesquare{228}{26}{28} parameter set $\mathcal{X}_1$ \hspace{2mm}
    \linecircle{55}{126}{184} parameter set $\mathcal{X}_2$ \hspace{2mm}
    \linedotted{0}{0}{0} time limit
    \vspace{-1mm}}
    \caption{Timing results for phase retrieval from \Subsec~\ref{subsec:phase_retrieval}.
    The verification problem for the final iterate of $\mathcal{X}_2$ reaches the time limit of $7200$ seconds (excluding time for OBBT).
    }
    \label{fig:phase_retrieval_timing_results}
\end{figure}

\begin{figure}[!h]
  \captionsetup{skip=4pt}
  \centering
  \includegraphics[width=0.44\linewidth]{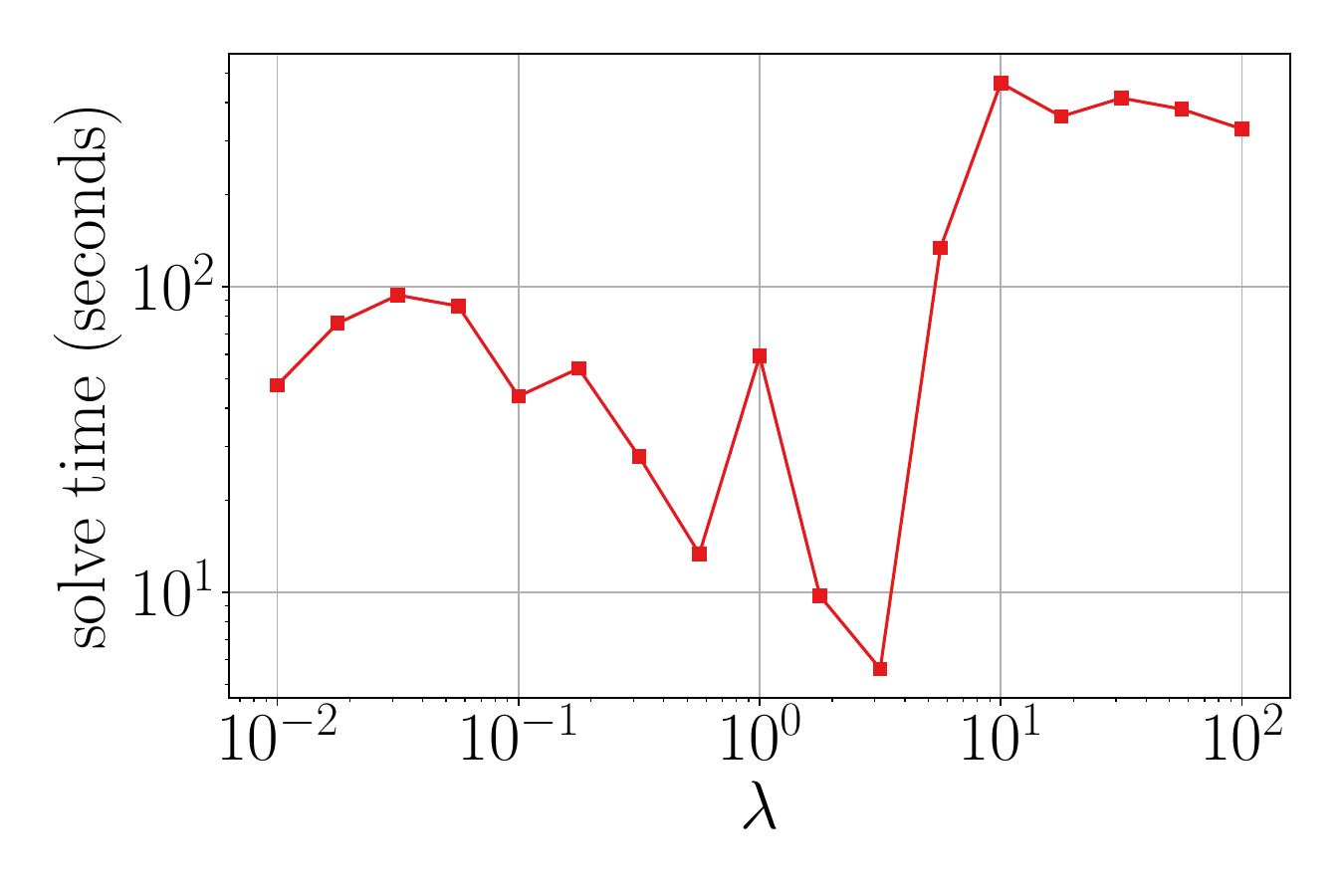}
  \vspace{-3mm}
    \caption{\vspace{-2mm}Timing results for sparse coding from \Subsec~\ref{subsec:sparse_coding}.
    }
    \label{fig:sparse_coding_timing_results}
\end{figure}

\begin{figure}[!h]
  \captionsetup{skip=4pt}
  \centering
  \includegraphics[width=0.88\linewidth]{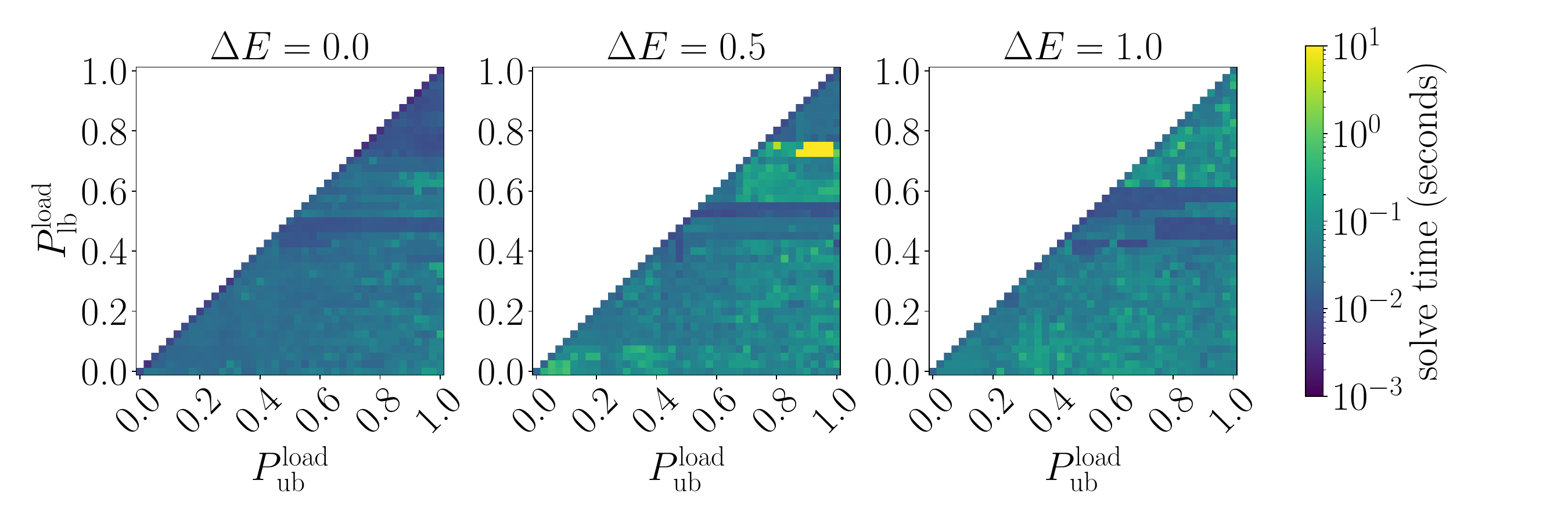}
  \vspace{-3mm}
    \caption{Timing results for hybrid vehicle control from \Subsec~\ref{subsec:hybrid_vehicle}.
    The time limit is $10$ seconds.
    }
    \label{fig:hybrid_vehicle_timing_results}
\end{figure}

We report solver times for the verification problems across all examples: box QP (Figure \ref{fig:box_qp_timing_results}), network utility (Figure \ref{fig:network_utility_timing_results}), power converter control (Figure \ref{fig:power_converter_timing_results}), knapsack (Figure \ref{fig:knapsack_timing_results}), phase retrieval (Figure \ref{fig:phase_retrieval_timing_results}), sparse coding (Figure \ref{fig:sparse_coding_timing_results}), and hybrid vehicle control (Figure~\ref{fig:hybrid_vehicle_timing_results}).

\end{document}